\documentclass[a4paper,11pt]{article}

\usepackage{etex}
\usepackage{amsmath}
\usepackage{amsfonts}
\usepackage{amssymb}
\usepackage{amsthm}
\usepackage{color}
\usepackage{subfig}
\usepackage[utf8]{inputenc}

\usepackage{pgfplots}
\usepackage{tikz}

\allowdisplaybreaks

\newtheorem{proposition}{Proposition}[section]
\newtheorem{theorem}[proposition]{Theorem}
\newtheorem{definition}[proposition]{Definition}

\theoremstyle{remark}

\definecolor{myBlue}{RGB}{0,131,204}
\definecolor{myGreen}{RGB}{153,192,0}
\definecolor{myRed}{RGB}{230,0,26}
\definecolor{myMagenta}{RGB}{0166,0,132}
\definecolor{myLilac}{RGB}{114,16,133}
\definecolor{myGrey}{RGB}{181,181,181}

\newcommand{\RR}{\mathbb R}
\newcommand{\N}{\mathbb{N}}
\newcommand{\R}{\mathbb{R}}
\newcommand{\vn}[1]{\| #1 \|}
\newcommand{\ST}{\quad{\text{s.t.}}\quad}

\newcommand{\expectedvalue}{\mu} 

\DeclareMathOperator*{\supp}{supp}

\DeclareMathOperator{\cov}{\textit{Cov}}
\DeclareMathOperator{\var}{VaR}
\DeclareMathOperator{\cvar}{CVaR}
\DeclareMathOperator{\rvar}{RVaR}
\DeclareMathOperator{\rcvar}{RCVaR}
\newcommand{\st}{\text{s.t. }}
\DeclareMathOperator*{\argmin}{argmin}

\begin{document}

\thispagestyle{empty}
\setcounter{page}{0}

\bigskip

\begin{center}
   {\bf Convergence of a Scholtes-type Regularization Method for Cardinality-Constrained Optimization Problems with an Application in Sparse Robust Portfolio Optimization}
\end{center}

\bigskip

\begin{center}
Martin Branda$^{1,2}$, Max Bucher$^3$, Michal \v{C}ervinka$^{1,4}$ and Alexandra Schwartz$^3$
\end{center}

\medskip

\noindent \hskip4.2cm $^1$Academy of Sciences of the Czech Republic

\noindent \hskip4.2cm Institute of Information Theory and Automation

\noindent \hskip4.2cm 
 Pod vod\'{a}renskou v\v{e}\v{z}\'{i} 4

\noindent \hskip4.2cm CZ-182 08 Prague 8 , Czech Republic

\noindent \hskip4.2cm e-mail: \{branda, cervinka\}@utia.cas.cz

\bigskip

\noindent \hskip4.2cm $^2$Charles University

\noindent \hskip4.2cm Faculty of Mathematics and Physics

\noindent \hskip4.2cm Department of Probability and Mathematical Statistics

\noindent \hskip4.2cm Sokolovsk\'{a} 83

\noindent \hskip4.2cm CZ-186 75 Prague 8, Czech Republic

\bigskip

\noindent \hskip4.2cm $^3$Technische Universität Darmstadt

\noindent \hskip4.2cm Graduate School Computational Engineering

\noindent \hskip4.2cm Dolivostra{\ss}e 15

\noindent \hskip4.2cm 64293 Darmstadt, Germany

\noindent \hskip4.2cm e-mail: \{bucher, schwartz\}@gsc.tu-darmstadt.de

\bigskip

\noindent \hskip4.2cm $^4$Charles University

\noindent \hskip4.2cm Faculty of Social Sciences

\noindent \hskip4.2cm Smetanovo n\'{a}b\v{r}e\v{z}\'i 6

\noindent \hskip4.2cm CZ-110 01 Prague 1, Czech Republic

\bigskip

\begin{center}
    March 30, 2017
\end{center}

\vfill

\begin{abstract}
We consider general nonlinear programming problems with cardinality constraints.
By relaxing the binary variables which appear in the natural mixed-integer programming formulation, we obtain an almost equivalent  nonlinear programming problem, which is thus still difficult to solve.
Therefore, we apply a Scholtes-type regularization method to obtain a sequence of easier to solve problems and investigate the convergence of the obtained KKT points.
We show that such a sequence converges to an S-stationary point, which corresponds to a local minimizer of the original problem under the assumption of convexity.

Additionally, we consider portfolio optimization problems where we minimize a risk measure under a cardinality constraint on the portfolio.
Various risk measures are considered, in particular Value-at-Risk and Conditional Value-at-Risk under normal distribution of returns and their robust counterparts under moment conditions.
For these investment problems formulated as nonlinear programming problems with cardinality constraints we perform a numerical study on a large number of simulated instances taken from the literature and illuminate the computational performance of the Scholtes-type regularization method in comparison to other considered solution approaches: a mixed-integer solver, a direct continuous reformulation solver and the Kanzow-Schwartz regularization method, which has already been applied to Markowitz portfolio problems.
\end{abstract}



\newpage

\section{Introduction}

Sparse solutions of optimization problems, i.e. solutions with a limited number of nonzero elements, are required in many areas including image and signal processing, mathematical statistics, machine learning, inverse environmental modeling and others.
The approach proposed in this paper is motivated mainly by applications of these problems in the portfolio selection theory where sparse solutions are popular due to lower fluctuations of the future (out-of-sample) portfolio performance, cf. DeMiguel et al. \cite{dMGNU}, and also due to the reduction of the transaction costs which are growing with the increasing number of assets included into the portfolio.
A standard way how to ensure the sparsity of the solutions is imposing a cardinality constraint where the number of nonzero elements of the solution is bounded.
While some studies employ a penalty on the $l_1$-norm of the asset weight vector or its alternatives, e.g. Fastrich et al. \cite{FPW}, some consider the explicit cardinality constraints.
The portfolio optimization problem resulting from the latter can be viewed as a mixed-integer problem and it is considered computationally challenging.
The examples of solution techniques include exact branch-and-bound methods, e.g., Borchers and Mitchell \cite{BM}, Bertsimas and Shioda \cite{BS}; exact branch-and-cut methods, e.g., Bienstock \cite{B}; heuristic algorithms, e.g., Chang et al. \cite{ChMBS}; and relaxation algorithms, e.g., Shaw et al. \cite{SLK}, Murray and Shek \cite{MS} and Burdakov et al. \cite{BKS2, BKS1}.
Portfolio optimization problems with sparsity have been one of the main motivations to study optimization problems with cardinality constraints, see \cite{DLR 12,RGS 10,SZL 13,ZSLS 13,BeE 13,FMPSW 13} and the references therein for some more ideas.

In this paper, we follow the approach from \cite{BKS1,CKS,FMPSW 13} and reformulate the cardinality constrained problem into a continuous optimization problem with orthogonality constraints.
The resulting structure is similar to mathematical programs with complementarity constraints (MPCC), see \cite{LPR 96,OKZ 98} and the references therein for some background.
Thus one can try to adapt solution methods originally designed for MPCCs to the continuous reformulation of cardinality constrained optimization problems.
One popular class of algorithms for MPCCs are the so-called regularization or relaxation methods, see for example \cite{Sch,StU 10,KS,LiF 05,KDB 09,DFN 05}.
In the recent paper \cite{BKS1} the regularization method from \cite{KS} was adapted to cardinality constrained problems. In this paper we want to focus on the older regularization method from \cite{Sch}, which was suggested for MPCCs by Scholtes in 2001.
Although the theoretical properties of this method are weaker than those of most of the other regularization methods, a numerical study on a collection of MPCC test problems in \cite{HKS 13} showed that the Scholtes regularization performs very well in practice.
Hence, in the theoretical part of this paper, we analyze the convergence properties of an adaptation of this regularization method to problems with cardinality constraints.
As it turns out, for cardinality constrained problems this regularization has competitive convergence properties compared for example to the method from \cite{BKS1}.
Furthermore we show that the regularized problems possess better properties than the original problem in terms of constraint qualifications.

In portfolio optimization, two basic types of decision-making frameworks are adopted: the utility maximization and the return-risk trade-off analysis, see, e.g., Levy \cite{Levy} for properties and relations between these two approaches. In the latter, it is important to define a risk that the concerned system has.
In optimization problems governed by uncertain inputs such as rates of return, typically represented as random variables, the risk is explicitly quantified by a risk measure.
In return-risk analysis, widely used both in theory and practice, an investor faces a trade-off between expected return and associated risk.
In his pioneering work in 1952, Markowitz \cite{M1} adopted variance as a measure of risk in his mean-variance analysis.
Many other alternatives were introduced since then.
Nowadays, Value-at-Risk (VaR), which measures the maximum loss that can be expected during a specific time horizon with a probability $\beta$ ($\beta$ close to 1), is widely used in the banking and insurance industry as a downside risk measure. Despite its popularity, VaR lacks some important mathematical properties.
Artzner at al. \cite{A} presented an axiomatic definition of risk measures and coined a coherent risk measure which has certain reasonable properties.
Conditional Value-at-Risk (CVaR), the mean value of losses that exceed the value of VaR, exhibits favorable mathematical properties such as coherence implying convexity.
Rockafellar and Uryasev \cite{RU, RU2} proposed to minimize CVaR for optimizing a portfolio so as to reduce the risk of high losses without prior computation of the corresponding VaR while computing VaR as a by-product.
Their CVaR minimization formulation results usually in convex or even linear programs which proved attractive for financial optimization and risk management in practice due to their tractability for larger real life instances.
For each of these risk measures, one can formulate corresponding mean-risk portfolio optimization problems.

Regardless of the risk measure used, these models are strongly dependent on the underlying distribution and its parameters, which are typically unknown and have to be estimated, cf. Fabozzi et al. \cite{FHZ}.
Investors usually face the so-called estimation risk as they rely on a limited amount of data to estimate the input parameters.
Portfolios constructed using these estimators perform very poorly in terms of their out-of-sample mean and variance as the resulting portfolio weights fluctuate substantially over time, cf. e.g. Michaud \cite{M3} and Chopra and Ziemba \cite{ChZ}.
As some reformulations of mean-risk portfolio problems depend on the assumption of normality, poor performance can also be caused by deviations of the empirical distribution of returns from normality.
One can thus also consider the distribution ambiguity in the sense that no knowledge of the return distribution for risky assets is assumed while the mean and variance/covariance are assumed to be known.
For these reasons, we examine portfolio policies based on robust estimators.
Robust portfolio selection deals with eliminating the impacts of estimation risk and/or distribution ambiguity.
Goldfarb and Iyengar \cite{GI} studied the robust portfolio in the mean-variance framework.
Instead of the precise information on the mean and the covariance matrix of asset returns, they introduced special types of uncertainty, such as box uncertainty and ellipsoidal uncertainty.
They also considered the robust VaR portfolio selection problem by assuming a normal distribution.
Chen et al. \cite{ChHZ} minimized the worst-case CVaR risk measure over all distributions with fixed first and second moment information.
The reader is referred to El Ghaoui et al. \cite{eGOO} and Popescu \cite{P} for other studies on portfolio optimization with distributional robustness.
Pa\c{c} and P{\i}nar \cite{PP} extend Chen et al. \cite{ChHZ} to the case where a risk-free asset is also included and distributional robustness is complemented with ellipsoidal mean return ambiguity. Other choices of the ambiguity set for VaR and CVaR are considered e.g. by T\"{u}t\"{u}nc\"{u} and Koening \cite{TK}, Pflug and Wozabal \cite{PW}, Zhu and Fukushima \cite{ZF}, DeMiguel and Nogales \cite{dMN} and Delage and Ye \cite{DY}.
For survey of recent approaches to construct robust portfolios, we refer to Kim et al. \cite{KKF}.
Despite the vast literature on robust portfolio optimization and many works on sparse portfolio optimization, there are only few works that concern both sparse and robust portfolios, cf. e.g. Bertsimas and Takeda \cite{BT}.

As an application of the general problem class, we consider the cardinality constrained minimization of VaR and CVaR under normality of asset returns and their robust counterparts under distribution ambiguity.
We assume that both first and second order moments of the random returns are known.
The resulting problem can then be solved using the following four approaches: solve a mixed-integer reformulation using GUROBI, solve the continuous reformulation directly using SNOPT, apply the Scholtes regularization and the regularization from Burdakov et al. \cite{BKS1}.
We perform a numerical experiment based on randomly generated test examples from the literature to compare these four approaches.
A similar numerical study has been reported in Burdakov et al. \cite{BKS1} for a cardinality constrained (and non-robust) mean-variance model, where the objective function was convex quadratic and the standard constraints linear.
Here, we investigate the investment problems with VaR and CVaR as introduced above, which leads to a more complicated convex objective function.
In order to be able to solve the resulting problem with GUROBI, the objective function has to be reformulated using a second-order cone constraint.

The paper is organized as follows:
We start Section \ref{sec:background} with a brief background on the continuous reformulation of cardinality constraints and the related optimality conditions and constraint qualifications.
Section \ref{sec:scholtes} is then devoted to the adaptation of the Scholtes regularization method and the analysis of its properties.
In Section \ref{sec:application}, we introduce the risk measures and define investment problems with a condition on portfolio sparsity.
Section \ref{sec:numerics} finally provides an extensive numerical comparison of all four afore mentioned solution approaches.

A few words on notation:
By $e \in \RR^n$ we denote the vector with all components equal to one.
For two vectors $x,y \in \RR^n$ the vector $x \circ y \in \RR^n$ denotes the componentwise (Hadamard) product of $x$ and $y$.
For a vector $x \in \R^n$ the number of nonzero components is denoted by $\|x\|_0$, i.e.
\[
   \|x\|_0 := |\{i=1,\ldots,n \mid x_i \neq 0\}|.
\]

\section{Cardinality Constrained Problems and a Continuous Reformulation}\label{sec:background}

In this section we want to provide the necessary background on cardinality constrained optimization problems and the continuous reformulation used in \cite{BKS1,CKS,FMPSW 13}.
Since the continuous reformulation has similarities to MPCCs, we also have to introduce suitable optimality conditions and constraint qualifications.

Let us consider a general cardinality constrained optimization problem
\begin{equation}\label{eq:CCgeneral}
   \min_x f(x)  \ST g(x) \leq 0, \ h(x) = 0, \ \|x\|_0 \leq \kappa,
\end{equation}
where $f:\R^n \to \R$, $g:\R^n \to \R^m$, and $h:\R^n \to \R^p$ are assumed to be continuously differentiable.
To simplify notation, we use the index sets
\[
   I_g(x) := \{i \mid g_i(x) = 0\}
   \qquad \text{and} \qquad
   I_0(x) := \{i \mid x_i = 0\}
\]
in the following.

Before we introduce the continuous reformulation on which our analysis is based, let us mention an alternative mixed integer reformulation, which was used for example in \cite{B} and which we will use in our numerical comparison to solve the portfolio optimization problems with GUROBI.
In case lower and upper bounds $l \leq x \leq u$ on the variable $x$ are known, problem \eqref{eq:CCgeneral} can be reformulated using binary decision variables into
\begin{equation} \label{mip}
   \begin{aligned}
      \min_{x,z} \,\,\, & f(x) \\
      \st & g(x) \leq 0, \ h(x) = 0, \\
      & z \in \{0, 1\}^n,\\
      & l \circ z \leq x \leq u \circ z,\\
      & e^\top z \leq \kappa.
   \end{aligned}
\end{equation}
If $x_i$ is nonzero, then the corresponding $z_i$ must be equal to one and by the reformulated cardinality constraint $e^\top z \leq \kappa$ this can happen at most $\kappa$ times.
However, as even for simple instances of cardinality constrained problems Bienstock \cite{B} showed the problem to be NP-complete, solving \eqref{mip} even using specialized global solution techniques can be computationally very time demanding.

Thus, we instead consider the following continuous reformulation of \eqref{eq:CCgeneral} introduced in Burdakov et al. \cite{BKS1}:
\begin{equation} \label{rp}
   \begin{aligned}
      \min_{x,y} \,\,\, & f(x) \\
      \st & g(x) \leq 0, \ h(x) = 0, \\
      & 0 \leq y \leq e,\\
      & x \circ y =0,\\
      & e^\top y \geq n-\kappa.
   \end{aligned}
\end{equation}
Here, in contrast to the previous reformulation, whenever $x_i$ is nonzero, the corresponding $y_i$ has to be equal to zero.
Due to the reformulated cardinality constraint $e^\top y \geq n-\kappa$ this can again occur at most $\kappa$ times.
Note that this problem is closely related to a mathematical program with complementarity constraints (MPCC) due to the ``half-complementarity''constraints $y \geq 0, x \circ y = 0$.
In case the additional constraint $x \geq 0$ is present as in the examples considered above, problem \eqref{rp} in fact is an MPCC.
Consequently \eqref{rp} is a nonconvex problem, even when the original cardinality constrained problem (except for the cardinality constraint of course) was convex.
Thus, one can in general not expect to obtain global minima.
But if one is for example interested in obtaining local solutions or good starting points for a global method, this approach can be useful.

For MPCCs in addition to the KKT conditions (or strong stationarity) several weaker optimality conditions such as Mordukhovich- and Clarke-stationarity are known, see for example \cite{CKS} for precise definitions.
In contrast to this it was shown in \cite{BKS1,CKS} that M-stationarity and all weaker concepts coincide for the continuous reformulation \eqref{rp} and that strong and M-stationarity can be reduced to the following conditions.

\begin{definition}
   Let $ (x^*, y^*) $ be feasible for \eqref{rp}.
   Then $ (x^*, y^*) $ is called
   \begin{itemize}
      \item[(a)] {\em S-stationary} (S = strong) if there exist multipliers $ \lambda \in \R^m, \mu \in \R^p $, and $ \gamma \in \R^n $ such that the following conditions hold:
         \begin{eqnarray*}
            \nabla f(x^*) + \sum_{i=1}^m \lambda_i \nabla g_i(x^*) +
            \sum_{i=1}^p \mu_i \nabla h_i(x^*) + \sum_{i=1}^n \gamma_i e_i = 0,
            & & \\
            \lambda_i \ge 0,  \ \lambda_i g_i(x^*) = 0 \quad
            \forall i = 1, \ldots, m, \\
            \gamma_i = 0 \quad \forall i \text{ such that } y_i^*=0.
         \end{eqnarray*}
      \item[(b)] {\em M-stationary} (M = Mordukhovich) if there exist multipliers  $ \lambda \in \R^m, \mu \in \R^p $, and $ \gamma \in  \R^n $ such that the following conditions hold:
         \begin{eqnarray*}
            \nabla f(x^*) + \sum_{i=1}^m \lambda_i \nabla g_i(x^*) +
            \sum_{i=1}^p \mu_i \nabla h_i(x^*) + \sum_{i=1}^n \gamma_i e_i = 0,
            & & \\
            \lambda_i \ge 0, \ \lambda_i g_i(x^*) = 0 \quad
            \forall i = 1, \ldots, m, \\
            \gamma_i = 0 \quad \forall i \text{ such that } x_i^* \neq 0.
         \end{eqnarray*}
   \end{itemize}
\end{definition}

Note that S-stationarity is equivalent to the KKT conditions of \eqref{rp} and still depends on the artificial variable $y$.
In contrast, the M-stationarity condition is slightly weaker but independent from the artificial variable $y$.
Following the idea from \cite{BKS1,CKS} one can also define constraint qualifications for \eqref{rp} depending on the original variable $x$ only.
This leads for example to the following version of Mangasarian-Fromowitz constraint qualification:

\begin{definition}
   Let $(x^*,y^*)$ be feasible for \eqref{rp}.
   Then $(x^*,y^*)$ satisfies the \emph{Cardinality Constrained Mangasarian-Fromowitz Constraint Qualification (CC-MFCQ)} if the gradients
   \begin{equation*}
      \nabla g_i(x^*)\ (i\in I_g(x^*))\quad\text{and}\quad \nabla h_i(x^*)\ (i=1,\dots,p),\ e_i\ (i\in I_0(x^*))
   \end{equation*}
   are positively linearly independent, i.e. if one \emph{cannot} find multipliers $\lambda \geq 0$ and $\mu, \gamma$ such that $(\lambda, \mu, \gamma) \neq 0$ and
   \[
      \sum_{i \in I_g(x^*)} \lambda_i \nabla g_i(x^*) +
      \sum_{i=1}^p \mu_i \nabla h_i(x^*) + \sum_{i\in I_0(x^*)} \gamma_i e_i = 0.
   \]
\end{definition}

It was shown in \cite{CKS} that every local minimum of \eqref{rp}, where CC-MFCQ or a weaker CC constraint qualification holds, is an S-stationary point.
This results differs from what is known for general MPCCs, where S-stationarity of local minima can only be guaranteed under an MPCC analogue of the linear independence constraint qualification.
Under MPCC-MFCQ and all weaker MPCC constraint qualifications, local minima of an MPCC can only be guaranteed to be M-stationary.

In this paper we are interested in portfolio optimization as an application, see Section \ref{sec:application} for details.
The resulting optimization problems turn out to be convex, except for the cardinality constraint of course.
For this special class of cardinality constrained optimization problems, it is known that S-stationarity of a point $(x^*,y^*)$ implies that it is a local minimum, see \cite{CKS}:

\begin{theorem}
	Consider \eqref{rp}, where $f, g_i:\R^n \to \R$ are convex and $h:\R^n \to \R^p$ is linear.
	Then every S-stationary point $(x^*,y^*)$ is a local minimum of \eqref{rp}.
\end{theorem}

\section{Properties of a Scholtes-type Regularization Method}\label{sec:scholtes}

Currently, there exist many different solution approaches for MPCCs, see \cite{LPR 96,OKZ 98} for an introduction.
One popular class of algorithms for MPCCs are so-called regularization methods, see \cite{Sch,StU 10,KDB 09,LiF 05,DFN 05,KS}, where one replaces the original, difficult problem by a sequence of simpler nonlinear programs (NLP), whose feasible set shrinks to the original one in the limit.
In this section, we briefly introduce the regularization method by Burdakov et al. in \cite{BKS1} for cardinality constrained problems.
Then we adapt the regularization method from Scholtes \cite{Sch} to cardinality constrained problems and analyze its convergence properties as well as the regularized subproblems.

In \cite{BKS1} the idea from Kanzow and Schwartz \cite{KS} was adapted to cardinality constrained problems:
The corresponding regularized problem was obtained by replacing the constraints $y \geq 0, x \circ y = 0$ in \eqref{rp} by the inequalities
\begin{equation} \label{KSregul}
   \Phi^+(x,y; t) \leq 0, \Phi^-(x,y; t) \leq 0, y \geq 0,
\end{equation}
where $\Phi^+_i(x,y;t) = \varphi(x_i,y_i;t)$ and $\Phi^-_i(x,y;t) = \varphi(-x_i,y_i;t)$ with
\[
   \varphi(a,b;t) = \begin{cases}
      (a-t)(b-t) & \mbox{ if } a+b \geq 2t, \\
      -\frac{1}{2}\left[(a-t)^2 + (b-t)^2 \right] & \mbox{ if } a+b < 2t.
   \end{cases}
\]
It is not difficult to see that for $t \geq 0$ the inequality $\varphi(a,b;t) \leq 0$ is equivalent to $\min\{a,b\} \leq t$.
In the case one knows $x \geq 0$, as in the considered application in portfolio optimization, one can of course eliminate the constraint $\Phi^-(x,y; t) \leq 0$.

In this paper, we want to adapt the regularization method introduced by Scholtes \cite{Sch} for MPCCs to \eqref{rp}.
Although this regularization technique is one of the oldest and the theoretical results known for MPCCs are weaker than those known for example for the regularization from \cite{KS}, it is numerically still very successful for MPCCs, see \cite{HKS 13}.

\begin{figure}[h]
   \centering
   \subfloat[Complementarity constraints]{
      \begin{tikzpicture}
         \draw[->] (-1.75,0) -- (0,0) node[below left]{$0$} -- (1.75,0) node[right]{$x_i$};
         \draw[->] (0,-0.5) -- (0,1) node[left]{$1$} -- (0,1.5) node[left]{$y_i$};
         \draw[-, line width = 3pt, myBlue] (0,1) -- (0,0) ;
         \draw[-, line width = 3pt, myBlue] (-1.5,0)  -- (1.5,0);
      \end{tikzpicture}
   }
   \quad
   \subfloat[Kanzow-Schwartz regularization]{
      \begin{tikzpicture}
         \fill[myBlue!50] (-1.5,0.5) -- (-0.5,0.5) -- (-0.5,1) -- (0.5,1) -- (0.5,0.5) -- (1.5,0.5) -- (1.5,0) -- (-1.5,0) -- cycle;
         \draw[->] (-1.75,0) -- (0,0) node[below left]{$0$} -- (0.5,0) node[below]{$t$} -- (1.75,0) node[right]{$x_i$};
         \draw[-, line width = 2pt, myBlue] (-1.5,0) -- (0,0) -- (1.5,0);
         \draw[-, line width = 2pt, myBlue] (-1.5,0.5) -- (-0.5,0.5) -- (-0.5,1) -- (0.5,1) -- (0.5,0.5) -- (1.5,0.5);
         \draw[->] (0,-0.5) -- (0,1) node[left]{$1$} -- (0,1.5) node[left]{$y_i$};
      \end{tikzpicture}
   }
   \quad
   \subfloat[Scholtes regularization]{
      \begin{tikzpicture}
         \fill[domain=0.25:1.5, smooth, variable=\x, myBlue!50] plot ({\x},{0.25/\x}) -- (1.5,0) -- (0,0) -- (0,1) -- cycle;
         \fill[domain=-1.5:-0.25, smooth, variable=\x, myBlue!50] plot ({\x},{-0.25/\x}) -- (0,1) -- (0,0) -- (-1.5,0) -- cycle;
         \draw[->] (-1.75,0) -- (0,0) node[below left]{$0$} -- (0.5,0) node[below]{$\sqrt{t}$} -- (1.75,0) node[right]{$x_i$};
         \draw[-, line width = 2pt, myBlue] (-1.5,0) -- (1.5,0);
         \draw[domain=0.25:1.5, smooth, variable=\x, line width = 2pt, myBlue] plot ({\x},{0.25/\x});
         \draw[domain=-1.5:-0.25, smooth, variable=\x, line width = 2pt, myBlue] plot ({\x},{-0.25/\x});
         \draw[line width = 2pt, myBlue] (-0.25,1) -- (0.25,1);
         \draw[->] (0,-0.5) -- (0,1) node[left]{$1$} -- (0,1.5) node[left]{$y_i$};
      \end{tikzpicture}
   }
   \caption{Illustration of the constraints $0 \leq y_i \leq 1, \ x_iy_i = 0$ and the two regularizations}\label{fig:illustrationRegularizations}
\end{figure}
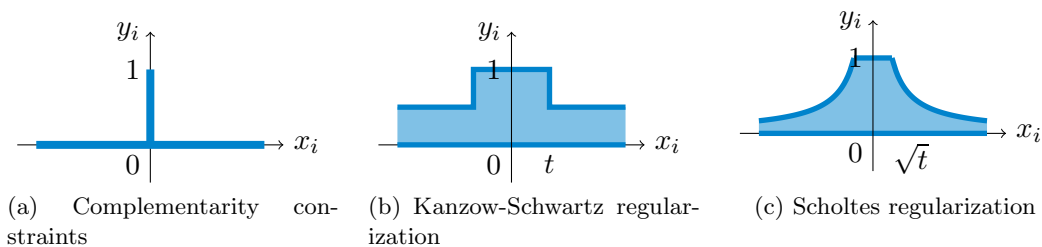

For a regularization parameter $t > 0$ we consider the regularized problem
\begin{eqnarray}\label{eq:nlpScholtes}
\begin{aligned}
   \text{NLP($t$)}:\quad\quad\min_{x,y} \,\,\, & f(x) \\
   \st & g(x) \leq 0, \quad h(x) = 0, \\
   & 0 \leq y \leq e,\\
   & -te \leq x \circ y \leq te,\\
   & e^\top y \geq n-\kappa.
\end{aligned}
\end{eqnarray}
Again it is easy to see that NLP($0$) corresponds to the original problem \eqref{rp} and that one can eliminate the constraint $x \circ y \geq -te$ in the case one knows $x \geq 0$.
A comparison of the feasible sets for on pair $(x_i,y_i)$ for both regularization techniques is given in Figure \ref{fig:illustrationRegularizations}.

Now consider a sequence of regularized problems NLP($t^k$) for $t^k \downarrow 0$ and assume that we can calculate a sequence of corresponding KKT points $(x^k,y^k)$, which is converging to some limit $(x^*,y^*)$.
Then one easily sees that the limit $(x^*,y^*)$ is feasible for the original problem \eqref{rp}.
However, the question remains whether $(x^*,y^*)$ is some kind of stationary point of \eqref{rp}, too.
In the classical MPCC case, one can only prove that the limit of such a sequence is a C-stationary point and there exist examples illustrating that this result is sharp.
Since C-stationary points and M-stationary points coincide for optimization problems with cardinality constraints, cf. \cite{CKS}, one would assume that we obtain M-stationarity of the limit here.
However, it turns out that the limit is in fact even S-stationary.

This result is even more surprising since it was shown in \cite{BKS1} that the adaptation of the Kanzow-Schwartz regularization retains its convergence properties knowns from MPCCs, i.e. converges to M-stationary points for both MPCCs and cardinality constrained problems.
But to be precise, the results for the Kanzow-Schwartz regularization were shown under a constant positive linear dependence constraint qualification, which is weaker than the MFCQ condition used for the Scholtes regularization.

\begin{theorem}
   Let $(t^k)_{k}$ be a sequence with $t^k>0$ for all $k\in\N$ and $t^k\downarrow0$.
   Let $(x^k,y^k)_{k}$ be a sequence of KKT points of \textnormal{NLP(}$t^k$\textnormal{)} converging to $(x^*,y^*)$.
   If CC-MFCQ holds at $(x^*,y^*)$, then $(x^*,y^*)$ is an S-stationary point of \eqref{rp}.
\end{theorem}

\begin{proof}
   Note that $(x^*,y^*)$ is a feasible feasible point of \eqref{rp}.
   Since $(x^k,y^k)_{k}$ is a sequence of KKT points of \textnormal{NLP(}$t^k$\textnormal{)}, there are multipliers $(\lambda^k,\mu^k,\tilde\gamma^k,\delta^k,\nu^k)$ for all $k\in\N$ such that
   \begin{subequations}
   \label{eq:KKT-grad}
   \begin{align}
      \nabla f(x^k)+\sum\limits_{i=1}^m\lambda_i^k \nabla g_i(x^k)+\sum\limits_{i=1}^p\mu_i^k\nabla h_i(x^k)+\sum\limits_{i=1}^n\tilde\gamma_i^ky_i^ke_i = 0, &\label{eq:kkt-grad-x}\\
      -\delta^ke+\sum\limits_{i=1}^n\nu_i^ke_i+\sum\limits_{i=1}^n\tilde\gamma_i^kx_i^ke_i = 0,&\label{eq:kkt-grad-y}\\
      \lambda_i^k\begin{cases}
         \geq 0, & \text{if } g_i(x^k)=0,\\
         = 0, & \text{otherwise},
      \end{cases} & \quad\forall i=1,\dots,m,\label{eq:kkt-relaxed-lambda}\\
      \delta^k\begin{cases}
         \geq 0, & \text{if } e^\top y^k = n-\kappa,\\
         = 0, & \text{otherwise},
      \end{cases}&\label{eq:kkt-relaxed-delta}\\
      \tilde\gamma_i^k\begin{cases}
         \geq0,&\text{if } x_i^k\cdot y_i^k=t^k,\\
         \leq0,&\text{if }x_i^k\cdot y_i^k=-t^k,\\
         =0,&\text{otherwise,}
         \end{cases}& \quad\forall i=1,\dots,n,\label{eq:kkt-relaxed-compl}\\
      \nu_i^k\begin{cases}
         \leq0,&\text{if } y_i^k=0,\\
         \geq0,&\text{if }y_i^k=1,\\
         =0,&\text{otherwise,}
         \end{cases}&\quad\forall i=1,\dots,n.\label{eq:kkt-bounds-y}
   \end{align}
   \end{subequations}
  Let us first have a closer look at the KKT conditions \eqref{eq:KKT-grad}.
  A componentwise inspection of equation \eqref{eq:kkt-grad-y} yields
   \begin{equation*}
      \delta^k=\nu_i^k+\tilde\gamma_i^kx_i^k
   \end{equation*}
   for all $i=1,\dots,n$.
   The sign restrictions on $\tilde\gamma^k$ imply $\tilde\gamma_i^k\cdot x_i^k\geq0$.
   Assuming there is an index $i\in\{1,\dots,n\}$ with $\nu_i^k<0$, it follows that $y_i^k=0$ and then using \eqref{eq:kkt-relaxed-compl} also $\tilde\gamma_i^k=0$.
   Thus the above equation yields $0>\nu_i^k=\delta^k\geq0$ which is a contradiction.
   Consequently we have
   \begin{equation}\label{eq:multipl-bounds-y-nonneg}
      \nu_i^k\geq0\quad\forall i=1,\dots,n.
   \end{equation}
   In case $\delta^k>0$ we have
   \begin{equation}\label{eq:aux5}
      0<\delta^k=\nu_i^k+\tilde\gamma_i^k x_i^k
   \end{equation}
   for all $i=1,\dots,n$.
   Then $\nu_i^k>0$ or $\tilde\gamma_i^kx_i^k>0$ has to hold for all $i=1,\dots,n$, which is true if and only if
   \begin{equation}\label{eq:aux1}
      y_i^k=1\quad\text{or}\quad y_i^k=\frac{t^k}{|x_i^k|}.
   \end{equation}
   For all $k \in \N$ define
   \[
      \gamma_i^k:= \tilde\gamma_i^ky_i^k  \quad \forall i=1,\dots,n.
   \]

   \smallskip\noindent
 \emph{Boundedness of the multipliers $(\lambda^k,\mu^k,\gamma^k)_{k}$:}
   We show this by contradiction.
   Thus, assume that $\lim_{k\rightarrow\infty} \vn{(\lambda^k,\mu^k,\gamma^k)} =\infty$. Then the sequence
   \begin{equation*}
      \left(\frac{(\lambda^k,\mu^k,\gamma^k)}{\vn{(\lambda^k,\mu^k,\gamma^k)}}\right)_{k\in\N}
   \end{equation*}
   is bounded and without loss of generality let it converge to some limit
   \begin{equation*}
      0\not=(\bar\lambda,\bar\mu,\bar\gamma):=\lim\limits_{k\rightarrow\infty}\frac{(\lambda^k,\mu^k,\gamma^k)}{\vn{(\lambda^k,\mu^k,\gamma^k)}}.
   \end{equation*}
   Clearly, $\bar\lambda\geq0$. Further, for all $i$ such that $g_i(x^*) < 0$ we have $g_i(x^k) < 0$ and thus also $\lambda_i^k = 0$ for all $k$ sufficiently large. That is, we have $\supp(\bar\lambda)\subseteq I_g(x^*)$.

   Next, to proceed with obtaining a contradiction, let us show that $\supp(\bar \gamma) \subseteq I_0(x^*)$.
   Assume, to the contrary, that there is an index $j\in\{1,\dots,n\}$, such that $x_j^*\not=0$ and $\bar\gamma_j\not=0$.
   Then we have $y_j^*=0$ and consequently
   \begin{equation*}
      x_j^k\not=0,\quad y_j^k < 1
   \end{equation*}
   for sufficiently large $k$.
   Since $\bar\gamma_j^k\not=0$ we have $\gamma_j^k\not=0$ and hence $\tilde\gamma_j^k\not=0$ for all $k$ sufficiently large.
   This implies $\delta^k=\nu_j^k+\tilde\gamma_j^kx_j^k>0$ and thus $\delta^k=\nu_i^k+\tilde\gamma_i^kx_i^k > 0$ for all $i=1,\dots,n$.
   Due to the KKT conditions, $\delta^k > 0$ is only possible if
   \begin{equation}\label{eq:aux2}
      e^\top y^k=n-\kappa
   \end{equation}
   for sufficiently large $k$.
   Furthermore, for sufficiently large $k$, $\tilde\gamma^k_j \neq 0$ implies
   \begin{equation}\label{eq:aux3}
     0 < y_j^k=\frac{t^k}{|x^k_j|} \,\, \text{ and } \,\, \nu_j^k=0
   \end{equation}

   As $y_j^k\rightarrow y_j^*=0$ and $y_j^k > 0$ hold for $k$ sufficiently large, the sequence $(y_j^k)_k$ is strictly monotonically decreasing (at least on a suitable subsequence).
   Moreover, since $e^\top y^k=n-\kappa$ for all $k$ sufficiently large (and $y^k$ is a finite-dimensional vector), strict monotone decrease of $(y_j^k)_k$ implies the existence of an index $l$ such that $(y_l^k)_k$ is strictly monotonically increasing (possibly on a suitable subsequence).
   Thus, we have
   \begin{equation*}
      y_l^*> 0,\ x_l^* = 0 \quad\text{and}\quad y_l^k<1,\ \nu_l^k=0\ \text{for sufficiently large }k.
   \end{equation*}
   Invoking \eqref{eq:aux1} and $y_l^k<1$ for sufficiently large $k$, we have
   \begin{equation}\label{eq:aux4}
      y_l^k=\frac{t^k}{|x_l^k|}.
   \end{equation}
   Since $\nu_j^k=\nu_l^k=0$, \eqref{eq:aux5}, \eqref{eq:aux3} and \eqref{eq:aux4} implies
   \begin{equation*}
      \frac{|\gamma_j^k|}{|\gamma_l^k|}=\frac{|\tilde\gamma_j^k\cdot y_j^k|}{|\tilde\gamma_l^k\cdot y_l^k|}
      =\frac{\left |\frac{\delta^k}{x_j^k}\cdot\frac{t^k}{|x_j^k|}\right |}{\left |\frac{\delta^k}{x_l^k}\cdot\frac{t^k}{|x_l^k|}\right |}
      =\left(\frac{x_l^k}{x_j^k}\right)^2 \xrightarrow[k\rightarrow\infty]{} \left(\frac{x_l^*}{x_j^*}\right)^2 = 0.
   \end{equation*}
   This leads to the contradiction
   \begin{equation*}
      0\not=|\bar\gamma_j|=\lim\limits_{k\rightarrow\infty}\frac{|\gamma_j^k|}{\vn{(\lambda^k,\mu^k,\gamma^k)}}
      \leq\lim\limits_{k\rightarrow\infty}\frac{|\gamma_j^k|}{|\gamma_l^k|}=0,
   \end{equation*}
   which concludes the proof of $\supp(\bar\gamma)\subseteq I_0(x^*)$.

   Now, dividing \eqref{eq:kkt-grad-x} by $\vn{(\lambda^k,\mu^k,\gamma^k)}$ and taking the limit $k\rightarrow\infty$ yields
   \begin{equation*}
      \sum\limits_{i\in I_g(x^*)}\bar \lambda_i \nabla g_i(x^*)+\sum\limits_{i=1}^p\bar\mu_i\nabla h_i(x^*)+\sum\limits_{i\in I_0(x^*)}\bar\gamma_ie_i = 0.
   \end{equation*}
   However, this, together with $\bar{\lambda}\geq0$, $\supp(\bar \lambda) \subseteq I_g(x^*)$, $\supp(\bar \gamma) \subseteq I_0(x^*)$, and $(\bar\lambda,\bar\mu,\bar\gamma)\not=0$, is in contradiction with the assumption of CC-MFCQ at $(x^*, y^*)$. Thus, the sequence of multipliers $(\lambda^k,\mu^k,\gamma^k)_{k}$ is bounded and without loss of generality we can assume that the whole sequence $(\lambda^k,\mu^k,\gamma^k)_{k}$ converges to some limit
   \begin{equation*}
      (\lambda^*,\mu^*,\gamma^*):=\lim\limits_{k\rightarrow\infty}(\lambda^k,\mu^k,\gamma^k).
   \end{equation*}
   Taking the limit in \eqref{eq:kkt-grad-x} as $k\rightarrow\infty$, we obtain
   \begin{equation*}
      \nabla f(x^*)+\sum\limits_{i\in I_g(x^*)}\lambda^*_i \nabla g_i(x^*)+\sum\limits_{i=1}^p\mu^*_i\nabla h_i(x^*)+\sum\limits_{i=1}^n\gamma^*_ie_i = 0.
   \end{equation*}

   \smallskip\noindent
   \emph{S-Stationarity of $x^*$ together with the multipliers $\lambda^*,\mu^*,\gamma^*$:}
   Using analogous arguments as previously, we have $\lambda^* \geq 0$ and $\supp(\lambda^*) \subseteq I_g(x^*)$.

   To complete the proof, it remains to show that $y_i^* = 0$ implies $ \gamma^*_i = 0$.
   Again, to the contrary, assume that there exists index $j$ such that $y_j^* = 0$ and $ \gamma^*_j \neq0$.
   This implies $ \gamma^k_j = \tilde\gamma^k_j y^k_j \neq 0$ for all $k$  sufficiently large and sequence $(y^k_j)_k$,
   \[
     0 < y^k_j = \frac{t^k}{|x^k_j|},
   \]
   is strictly monotonically decreasing to zero (at least on a suitable subsequence).
   Thus, we have $x_j^k \neq 0$ and $\nu^k_j = 0$ for all $k$ sufficiently large which together with $\gamma^k_j  \neq 0$  implies $\delta^k = \tilde\gamma^k_jx^k_j > 0$ and $e^\top y^k = n-\kappa$ for all $k$  sufficiently large. Analogously to the previous part of the proof, there has to exist an index $l$ such that $(y^k_l)_k$ is strictly increasing and
   \[
      0< y^k_l = \frac{t^k}{|x^k_l|}.
   \]
   This implies $y^k_l \to y^*_l > 0$ and $\nu^k_l = 0$ and thus $\delta^k = \gamma^k_l x^k_l$ and $x^k_l \neq 0$ for all $k$ sufficiently large.
   Finally, this together with $\gamma_i^k = \tilde\gamma^k_i y^k_i$ for all $i$ yields
   \begin{equation*}
      \frac{|\gamma^*_j|}{|\gamma^*_l|}
      =  \lim\limits_{k\rightarrow\infty}\frac{|\gamma_j^k|}{|\gamma_l^k|}
      =  \lim\limits_{k\rightarrow\infty}\frac{|\tilde\gamma_j^ky_j^k|}{|\tilde\gamma_l^ky_l^k|}
      =  \lim\limits_{k\rightarrow\infty}\frac{\left |\frac{\delta^k}{|x_j^k|}t^ky_j^k\right|}{\left |\frac{\delta^k}{|x_l^k|}t^ky_l^k\right|}
      = \lim\limits_{k\rightarrow\infty}\left(\frac{y_j^k}{y_l^k}\right)^2
      = \left(\frac{y_j^*}{y_l^*}\right)^2 = 0,
   \end{equation*}
   which is a contradiction to $ \gamma_j \neq 0$.
   This completes the proof.
\end{proof}

The previous proof differs slightly from the one typically used in the MPCC case, see e.g. \cite{HKS 13}, since CC-MFCQ does not demand positive linear independence of the gradients of the constraints with respect to $y$.
This is the reason why we only normalized the multipliers $\lambda^k, \mu^k, \gamma^k$ corresponding to constraints on $x$ in the previous proof.
The drawback of this approach is that the verification of the correct support of the limit of $\gamma^k$ is then more lengthy.
The central idea to exploit the fact $e^\top y^k = n- \kappa$ was borrowed from \cite{AB}, where a similar structure was used to reformulate chance constrained optimization problems.

Instead of normalizing only the multipliers $\lambda^k, \mu^k, \gamma^k$ one can also normalize all multipliers $\lambda^k, \mu^k, \gamma^k, \nu^k, \delta^k$.
This simplifies the verification of the correct support of $\gamma^k$ but makes it harder to obtain a contradiction to CC-MFCQ.
We follow this route in the proof of the next result, where we show that the regularized problems NLP($t$) indeed have better properties than the original problem \eqref{rp}.

\begin{theorem}
   Let $(x^*,y^*)$ be feasible for \eqref{rp} and CC-MFCQ hold there.
   Then there is a neighborhood $\mathcal{U}$ of $(x^*,y^*)$ such that for all $t > 0$ the standard MFCQ for NLP($t$) holds at every $(x,y)\in \mathcal{U}$ feasible for NLP($t$).
\end{theorem}

\begin{proof}
   Let us assume that the claim is false.
   Then there exist sequences $(x^k,y^k)_k \to (x^*,y^*)$ and $(t^k)_k > 0$ such that $(x^k,y^k)$ is feasible for NLP($t^k$) but MFCQ is violated.
   Consequently, we can find multipliers
$(\lambda^k,\mu^k,\tilde\gamma^k,\delta^k,\nu^k)_k$
   such that for all $k\in\N$
   \begin{equation}\label{eq:mfcq-coeff-not-zerovector}
      (\lambda^k,\mu^k,\tilde\gamma^k,\delta^k,\nu^k) \not=0\quad\forall k\in\N,
   \end{equation}
   along with condition
   \begin{align}
      \sum\limits_{i=1}^m\lambda_i^k \nabla g_i(x^k)+\sum\limits_{i=1}^p\mu_i^k\nabla h_i(x^k)+\sum\limits_{i=1}^n\tilde\gamma_i^ky_i^ke_i = 0, &\label{eq:mfcq-lincomb-grad-x}\\
   \end{align}
   and \eqref{eq:kkt-grad-y}-\eqref{eq:kkt-bounds-y} are satisfied.
   This means the relevant gradients are positively linearly dependent at $(x^k,y^k)$, hence MFCQ is violated.
Analogously to the proof of Theorem 3.1, we can show that $\nu_i^k < 0$ cannot occur, thus we have to consider $\nu_i^k \geq 0$ only.

   Now, let us define $\gamma_i^k:=\tilde\gamma_i^ky_i^k$ for all $i=1,\dots,n$ and all $k\in\N$.
   Because for all $i=1,\dots,n$ and all $k\in\N$ we have $\tilde\gamma_i^k\not=0$ only if $y_i^k\not=0$, it follows that
   \begin{equation*}
      \forall i=1,\dots,n,\ \forall k\in\N:\quad\left[\gamma_i^k\not=0\quad\Leftrightarrow\quad\tilde\gamma_i^k\not=0\right].
   \end{equation*}
   Using \eqref{eq:kkt-relaxed-compl} we have
   \begin{equation*}
      \tilde\gamma_i^kx_i^k=
      \begin{cases}
         \gamma_i^k\cdot\frac{x_i^k}{y_i^k},&\text{if } x_i^k\cdot y_i^k=|t|,\\
         0,&\text{otherwise.}
      \end{cases}
   \end{equation*}
   Therefore,
 we can write equations \eqref{eq:mfcq-lincomb-grad-x} and \eqref{eq:kkt-grad-y} for all $k\in\N$ as
   \begin{align}
      \sum\limits_{i=1}^m\lambda_i^k \nabla g_i(x^k)+\sum\limits_{i=1}^p\mu_i^k\nabla h_i(x^k)+\sum\limits_{i=1}^n\gamma_i^ke_i = 0, &\label{eq:mfcq-NEW-lincomb-grad-x}\\
      -\delta^ke+\sum\limits_{i=1}^n\nu_i^ke_i+\sum_{x_i^k\cdot y_i^k=|t|}\gamma_i^k\frac{x_i^k}{y_i^k}e_i = 0.&\label{eq:mfcq-NEW-lincomb-grad-y}
   \end{align}
   Due to the sign restrictions in \eqref{eq:kkt-relaxed-compl} we know $\tilde\gamma_i^kx_i^k\geq0$ for all $i=1,\dots,n$ and all $k\in\N$.
   Hence we also have
   \begin{equation}\label{eq:mfcq-product-new-coeff-x-nonnegative}
      \gamma_i^kx_i^k=\tilde\gamma_i^kx_i^ky_i^k\geq0\quad\forall i=1,\dots,n,\ \forall k\in\N.
   \end{equation}

   We proceed to deduce a contradiction with CC-MFCQ at $(x^*,y^*)$.
   Since by assumption
 $(\lambda^k,\mu^k,\tilde\gamma^k,\delta^k,\nu^k) \not=0$ for all $k\in\N$,
 we can choose the multipliers without loss of generality such that $\vn{(\lambda^k,\mu^k,\gamma^k,\delta^k,\nu^k)}=1$ for all $k\in\N$ and that the whole sequence converges to some limit
   \begin{equation*}
      0\not=(\lambda,\mu,\gamma,\delta,\nu):=\lim\limits_{k\rightarrow\infty}(\lambda^k,\mu^k,\gamma^k,\delta^k,\nu^k).
   \end{equation*}

   We have $\lambda\geq0$.
   Since for all $i$ such that $g_i(x^*) < 0$ we know $g_i(x^k) < 0$ and thus $\lambda_i^k = 0$ for all $k$ sufficiently large, we have
   \begin{equation}\label{eq:mfcq-support-lambda}
      \supp(\lambda)\subseteq I_g(x^*).
   \end{equation}

   We will prove $\supp(\gamma)\subseteq I_0(x^*)$ by contradiction.
   To this end we assume that there is an index $j\in\{1,\dots,n\}$ such that $\gamma_j\not=0$ and $x_j^*\not=0$.
   This implies $x_j^k\cdot y_j^k=|t|$ for sufficiently large $k$ and $y_j^*=0$.
   Therefore we know $x_j^k\not=0$ and $y_j^k>0$ for $k$ sufficiently large and $y_j^k\rightarrow0$ for $k\rightarrow\infty$.
   Thus we have $y_j^k<1$ and hence $\nu_j^k=0$ for sufficiently large $k$.
   Keeping in mind \eqref{eq:mfcq-product-new-coeff-x-nonnegative} it follows from the $j$-th row of \eqref{eq:mfcq-NEW-lincomb-grad-y} that
   \begin{equation*}
      \delta^k=\nu_j^k+\gamma_j^k\frac{x_j^k}{y_j^k}=\gamma_j^k\frac{x_j^k}{y_j^k}\longrightarrow\infty\quad(k\rightarrow\infty).
   \end{equation*}
   Because $(\lambda^k,\mu^k,\gamma^k,\delta^k,\nu^k)_{k
   }$ is convergent, this is a contradiction.
   Consequently we have
   \begin{equation}\label{eq:mfcq-support-gamma}
      \supp(\gamma)\subseteq I_0(x^*).
   \end{equation}

   It remains to show that $(\lambda,\mu,\gamma)\not=0$.
   We will show this also by contradiction.
   Assume that $(\lambda,\mu,\gamma)=0$.
   Since $(\lambda,\mu,\gamma,\delta,\nu)\not=0$ this implies $(\delta,\nu)\not=0$.
   If $y_i^* =0 $, then we have $y_i^k<1$ and hence $\nu_i^k=0$ for sufficiently large $k$.
   Thus we have
   \begin{equation*}\label{eq:mfcq-aux1}
      \nu_i=0\quad\forall i\,:\,y_i^* = 0.
   \end{equation*}
   If $y_i^*\neq0$, we have $y_i^k > 0$ for sufficiently large $k$ and $x^*_i = 0$.
   Using \eqref{eq:mfcq-NEW-lincomb-grad-y} and the assumption $\gamma = 0$ it follows that
   \begin{equation*}
      \delta = \lim\limits_{k\rightarrow\infty} \delta^k= \lim\limits_{k\rightarrow\infty} \nu_i^k+\gamma_i^k\frac{x_i^k}{y_i^k}=\nu_i+\gamma_i\frac{x_i^*}{y_i^*}=\nu_i.
   \end{equation*}
   Because $(\delta,\nu)\not=0$ we consequently have
   \begin{equation}\label{eq:mfcq-aux2}
      \delta=\nu_i>0\quad\forall i\,:\,y_i^*\neq 0.
   \end{equation}
   This implies $\nu_i^k>0$ and thus $y_i^k=1$ for $k$ sufficiently large for all $i$ such that $y^*_i \neq 0$.
   Thus, we have
   \begin{equation}\label{eq:mfcq-aux3}
      y_i^*\in\{0,1\}\quad\forall i=1,\dots,n.
   \end{equation}
   It also follows from \eqref{eq:mfcq-aux2} that $\delta^k>0$ and hence for all $k$ sufficiently large we have
   \begin{equation}\label{eq:mfcq-aux4}
      e^\top y^k=n-\kappa.
   \end{equation}
   Thus $e^\top y^*=n-\kappa$ holds.
   Since \eqref{eq:mfcq-aux3} holds, we consequently have
   \begin{equation*}
      \left|\{i\,:\,y_i^*\not=0\}\right|=n-\kappa
      \qquad \text{and} \qquad
      \left|\{i\,:\,y_i^*=0\}\right|=\kappa.
   \end{equation*}
   This implies that for all $k$ sufficiently large
   \begin{equation*}
      n-\kappa = e^\top y^k
       = \sum\limits_{y_i^*\not=0}y_i^k+\sum\limits_{y_i^*=0}y_i^k
       = n-\kappa+\sum\limits_{y_i^*=0}y_i^k.
   \end{equation*}
   Because $y^k\geq0$, it follows from the above equation that $y_i^k=0$ for all $i$ such that $y_i^*=0$ and thus $y^k = y^*$ for all $k$ sufficiently large.

   For the $\kappa >0$ indices $i$ with $y^*_i = y^k_i = 0$ we have $\nu_i^k=\gamma_i^k=0$.
   Then \eqref{eq:mfcq-NEW-lincomb-grad-y} yields $\delta^k=0$, a contradiction to \eqref{eq:mfcq-aux2}.
   Thus the assumption $(\lambda,\mu,\gamma)=0$ is false and we have
   \begin{equation*}
      (\lambda,\mu,\gamma)\not=0.
   \end{equation*}

   Using \eqref{eq:mfcq-support-lambda} and \eqref{eq:mfcq-support-gamma}, it follows from \eqref{eq:mfcq-NEW-lincomb-grad-x} for $k\rightarrow\infty$ that
   \begin{equation}
      \sum\limits_{i\in I_g(x^*)}\lambda_i \nabla g_i(x^*)+\sum\limits_{i=1}^p\mu_i\nabla h_i(x^*)+\sum\limits_{i\in I_0(x^*)}\gamma_ie_i = 0.
   \end{equation}
   Since $(\lambda,\mu,\gamma)\not=0$ and $\lambda\geq0$, this is a contradiction to CC-MFCQ.
\end{proof}

\section{Robust Risk Measures for Portfolio Optimization under Distribution Ambiguity}\label{sec:application}

In this section we want to provide an application for the abstract cardinality constrained optimization problems discussed in the previous sections.
To do so, we consider the following portfolio optimization problem:
\begin{equation} \label{gp}
   \begin{aligned}
      \min_x \,\,\, & r(x) \\
      \st & e^\top x = 1, \\
      & 0 \leq x \leq u,\\
      & \|x\|_0 \leq \kappa.
   \end{aligned}
\end{equation}
Here, we consider a market with $n$ risky financial assets.
The disposable wealth is to be allocated into a portfolio $x \in \R^n$, such that each component $x_i$ denotes the fraction of disposable wealth to be invested into the $i$-th asset for $i=1, \dots, n$.
We do not allow short-sales, i.e. we assume $x \geq 0$.
For numerical purposes, we assume that there exist upper bounds $u \geq 0$ on the possible investment.
If no upper bounds are present, one can use $u = e$ due to the budget constraint $e^\top x = 1$.
The latter states that  we demand that the whole disposable budget is invested.
Additionally, we introduce the cardinality constraint $\|x\|_0 \leq \kappa$, i.e. one may invest in at most $\kappa$ assets.
Naturally, we assume that $\kappa < n$.

For a vector $x$ of allocations to $n$ risky assets and a random vector $\xi$ of return rates for these assets, we consider the following loss function
\[
   \ell(x,\xi) = - x^\top \xi.
\]
Assume that $\xi$ follows a probability distribution $\pi$ from the ambiguity (uncertainty) set
$D = \{\pi \mid E_\pi[\xi] = \expectedvalue, \cov_\pi(\xi) = \Gamma \succ 0 \}$ of distributions with expected value $\expectedvalue$ and positive semidefinite covariance matrix $\Gamma$. Note that Markowitz \cite{M1} considered the variance $\sigma^2(x) = x^\top \Gamma x$ as a risk measure associated with a portfolio $x$.

In the 90s, the investment bank J.P. Morgan reinvented the quantile risk measure (quantile premium principle) used by actuaries for investment banking, giving rise to Value-at-Risk (VaR).
Associated with a confidence level $\beta$ and portfolio $x$, VaR is defined as
\[
   \var_\beta(x) = \min \{z \mid P_\pi( \ell(x, \xi) \leq z) \geq \beta\}.
\]

Artzner et al. \cite{A} defined coherent risk measure as a measure satisfying monotonicity, translation invariance, subadditivity and positive homogeneity.
It is known, that VaR is not a coherent risk measure as it fails subadditivity.
On the other hand, the conditional value-at-risk (CVaR), as introduced by Rockafellar and Uryasev \cite{RU}, turns out to be a convex and coherent risk measure.
CVaR at level $\beta$ is defined as the expected value of loss that exceeds $\var_\beta(x)$.
Alternatively, Rockafellar and Uryasev \cite{RU} showed that calculation of CVaR and VaR can be achieved simultaneously by minimizing the following auxiliary function with respect to $\alpha \in \RR$
\[
   F_\beta(x,\alpha) = \alpha + \frac{1}{1-\beta}E[(\ell(x,\xi)-\alpha)_+],
\]
where $(v)_+ = \max\{0, v\}$.
Thus,
\[
   \cvar_\beta(x) = \min_\alpha F_\beta(x, \alpha)
\]
and $\var_\beta(x)$ is the left endpoint of the interval $\argmin_\alpha F_\beta(x, \alpha)$.

Let us assume normality of returns $\xi$.
Denote by $\phi$ and $\Phi$ density and cumulative distribution function of the standard normal distribution, respectively.
Following Fabozzi et al. \cite{FHZ}, originating in Rockafellar and Uryasev \cite{RU}, the Value-at-Risk can be expressed as
\begin{equation}\label{VaR}
   \var_\beta(x) = \zeta_\beta \sqrt{x^\top Q x} - \expectedvalue^\top x,
\end{equation}
where $\zeta_\beta = -\Phi^{-1}(1-\beta)$, and assuming $\beta > 0.5$, the Conditional Value-at-Risk reduces to
\begin{equation}\label{CVaR}
   \cvar_\beta(x) = \eta_\beta \sqrt{x^\top Q x} - \expectedvalue^\top x,
\end{equation}
where $\eta_\beta = \frac{-\int_{-\infty}^{\Phi^{-1}(1-\beta)}t\phi(t)dt}{1-\beta}$.

Further, we consider the worst case VaR for a fixed $x$ with respect to the ambiguity set $D$ defined as
\[
   \rvar_\beta(x) = \sup_{\pi \in D} \var_\beta(x).
\]
Analogously, we consider the worst case CVaR for a fixed $x$ with respect to set $D$ defined as
\[
   \rcvar_\beta(x) = \sup_{\pi \in D} \cvar_\beta(x) = \sup_{\pi \in D} \min_\alpha F_\beta(x, \alpha).
\]

Based on Chen et. al \cite[proof of Theorem 2.9]{ChHZ}, further generalized in Pa\c{c} and P{\i}nar \cite{PP} using Shapiro \cite[Theorem 2.4]{S}, we obtain the following representations that under distribution ambiguity,
\begin{equation}\label{RVaR}
   \rvar_\beta(x) = \frac{2\beta-1}{2\sqrt{\beta(1-\beta)}} \sqrt{x^\top Q x} - \expectedvalue^\top x
\end{equation}
and
\begin{equation}\label{RCVaR}
   \rcvar_\beta(x) = \sqrt{\frac{\beta}{1-\beta}} \sqrt{x^\top Q x} - \expectedvalue^\top x.
\end{equation}

In the following section we consider cardinality constrained portfolio selection models for each of the risk measures \eqref{VaR}--\eqref{RCVaR} replacing the general risk function $r(x)$ in \eqref{gp}.

\section{Numerical Comparison of Different Solution Methods}\label{sec:numerics}

In this section, we compare the performance of the Scholtes regularization method introduced in this paper for cardinality constrained optimization problems with the Kanzow-Schwartz regularization from \cite{BKS1}, the direct application of an NLP solver to the continuous reformulation \eqref{rp} and the solution of \eqref{mip} with a mixed integer solver.
We test all four approaches on the investment problems described in the previous section with the VaR and CVaR measures under normality assumption and the robust VaR and CVaR.
Moreover, we consider each problem for several levels of $\beta$, in particular we select $\beta \in \{0.9, 0.95, 0.99 \}$.
Table \ref{tab:0} contains the values of the corresponding quantiles and generalized quantiles, which appear in the exact reformulations of the risk measures.

\begin{table}
   \caption{Quantiles and generalized quantiles as defined in \eqref{VaR}--\eqref{RCVaR} }\label{tab:0}
   \centering
   \begin{tabular}{cc|c|c|c}
      & $c_\beta \diagdown \beta$ & 0.9 & 0.95 & 0.99\\\hline
      VaR & $\zeta_\beta$ & 1.2816 & 1.6449 & 2.3263\\
      CVaR & $\eta_\beta$& 1.7550 & 2.0627 & 2.6652\\
      RVaR & $\frac{2 \beta -1}{2\sqrt{\beta(1-\beta)}}$& 1.3333 & 2.0647 & 4.9247\\
      RCVaR & $\sqrt{\frac{\beta}{1-\beta}}$ & 3.0000 & 4.3589 & 9.9499\\
   \end{tabular}
\end{table}

We use 90 simulated instances with mean vectors and variance matrices, which were already employed by \cite{FG} and are freely available at website \cite{FGwp}.
The generation of the data was described by \cite{Pard}.
For each number $n = 200, 300$ and $400$ of assets there are $30$ different problems included in the dataset.

We compare the performance of the following solution approaches:
\begin{enumerate}
   \itemsep -3pt
	\item \texttt{GUROBI\_60}: Solve the mixed integer formulation \eqref{mip} using the commercial mixed-integer solver \texttt{GUROBI}, version 6.5, with time limit of $60$s and start vector $x^0 = 0$, $z^0 = e$.

	\item \texttt{GUROBI\_300\_40}: Same as above but with time limit $300$s and node limit $40$.

	\item \texttt{Relaxation\_01}: Solve the continuous reformulation \eqref{rp} directly using the sparse SQP method \texttt{SNOPT}, version 7.5, with start vectors $x^0 = 0$, $y^0 = e$.

	\item \texttt{Relaxation\_00}: Same as above but with start vectors $x^0 = 0$, $y^0 = 0$.

	\item \texttt{Scholtes\_01}: Solve a sequence of Scholtes regularizations \eqref{eq:nlpScholtes} using \texttt{SNOPT} with starting point $x^0 = 0$, $y^0 = e$.

	\item \texttt{Scholtes\_00}: Same as above but with start vectors $x^0 = 0$, $y^0 = 0$.

   \item \texttt{KanzowSchwartz\_01}: Solve a sequence of Kanzow--Schwartz regularizations \eqref{KSregul} using \texttt{SNOPT} with starting point $x^0 = 0$, $y^0 = e$.

   \item \texttt{KanzowSchwartz\_00}: Same as above but with start vectors $x^0 = 0$, $y^0 = 0$.
\end{enumerate}

All computations were done in \texttt{MATLAB R2014a}.
Before we discuss the results, let us state a few details on the implementation of the respective solution approaches:

More information on the solver \texttt{GUROBI} and its various options can be found at \cite{G}.
To be able to solve the mixed-integer problem \eqref{mip} with \texttt{GUROBI}, we had to reformulate it in the following form:
\begin{equation*}
   \begin{aligned}
      \min_{x,z,w,v} \,\,\, &  c_\beta v  - \expectedvalue^\top x \\
      \st & e^\top x = 1,\\
      & 0 \leq x \leq u \circ z,\\
      & z \in \{0, 1\}^n,\\
      & e^\top z \leq \kappa,\\
      & v \geq 0,\\
      & w = Q^{\frac{1}{2}} x,\\
      & v^2 \geq w^\top w,
   \end{aligned}
\end{equation*}
where $c_\beta$ is the respective constant from Table \ref{tab:0} for the different risk measures.
Since we used $x^0 = 0$ as start vector, we also used $w^0 = 0$ and $v^0 = 0$.

Note that \texttt{GUROBI} is a global solver, i.e. it tries to verify that a candidate solution is indeed a global minimum.
Since the other solution approaches do not provide any guarantee of finding a global solution, we set the option \texttt{mipfocus} to 1 in order to encourage \texttt{GUROBI} to try to find good solutions fast.
Additionally we set the option \texttt{timelimit} to $60$s at first.
However, we observed that \texttt{GUROBI} sometimes spent the whole time by looking for a feasible solution without moving to the branch-and-bound tree.
Thus we increased \texttt{timelimit} to $300$s and added the condition on the maximal number of computed nodes $\texttt{nodelimit} = 40$ to obtain results less dependent on slight variations in computation time.

The continuous reformulation \eqref{rp} and the regularized problems \eqref{KSregul} and \eqref{eq:nlpScholtes} were all solved using the sparse SQP method \texttt{SNOPT}, see \cite{GMS, GMSW} for more information.
We started both regularization methods with an initial parameter $t^0 = 1$ and decreased $t^k$ in each iteration by a factor of $0.01$.
Both regularization methods were terminated if either $\|x^k \circ y^k\|_\infty \leq 10^{-6}$ or $t^k < 10^{-8}$.

The constraints $e^\top x = 0$ and $0 \leq x \leq u$ were usually satisfied in the solutions $x^*$ found by all methods (except for \texttt{GUROBI}, which occasionally did not return a feasible solution at all, see below).
In order to check whether the cardinality constraint $\|x\|_0 \leq \kappa$ are also satisfied, we counted the number of all components $x^*_i > 10^{-6}$.

Table \ref{tab:1} contains results for a particular problem with 400 assets (pard400-e-400).
We can see that \texttt{GUROBI} running 60s was not able to provide a feasible solution for problem with $\rvar_{0.99}$.
The Scholtes regularization starting from point $x^0 = 0$, $y^0= e$ was not successful for $\rcvar_{0.95}$.
However, in all other cases the Scholtes regularization starting from $x^0 = 0$, $y^0 = e$ provided the best solution with a runtime around $1s$.
In Table \ref{tab:1} we also report the relative gap
$$
   \text{relative  gap} = (f-f_\text{best})/f_\text{best},
$$
where $f$ is the objective value obtained by an algorithm and $f_\text{best}$ denotes the lowest objective value for a problem.

Summary results for all problems are reported in Tables \ref{tab:3}, \ref{tab:4}, \ref{tab:5}.
For each problem with a particular risk measure, level $\beta$, number of assets and algorithm we report the following descriptive statistics over 30 instances of problems:
average relative gap with respect to the minimal objective value, average computation time (in seconds), number of cases when the algorithm found the best solution, number of cases when the result was infeasible with respect to the sparsity or orthogonality.
All computations were done on two computers with comparable performance indicators.
Nonetheless, the given computation times should only be used for a qualitative comparison of the methods.

If we consider the initial value $x^0 = 0$, $y^0 = e$, it can be observed that the best results were obtained by the Scholtes regularization (Alg. 5).
When the results of this regularization were feasible, they correspond to the best obtained solutions.
However, in a few cases the portfolios obtained by the Scholtes regularization and in more cases the results obtained by the Kanzow-Schwartz regularization (Alg. 7) were infeasible.
Also \texttt{SNOPT} applied directly to the continuous reformulation (Alg. 3) behaved very badly showing an average relative gap greater than 100\%.

To further investigate the behavior, we changed the starting point of the continuous reformulation and both regularizations to $x^0 = 0$, $y^0 = 0$.
In this case, the obtained optimal values were slightly worse for both regularizations, but we have reduced the problems with infeasible solutions.
Moreover, for the starting point $x^0 = 0$, $y^0 = 0$, the behavior of the continuous reformulation approach (Alg. 4) improved significantly such that it is fully comparable with the regularizations.

Figures \ref{fig:1}, \ref{fig:2}, \ref{fig:3} present performance plots for each problem size and algorithm.
We identified the minimal objective value for each problem found by any of the eight considered algorithms and then compared it with the remaining objective values using the ratio: actual objective value/minimal objective value.
The graphs report the relative number of problems ($y$-value), where the ratio is lower or equal to the $x$-value.
We would prefer algorithms with the curve close to the upper-left corner, i.e. which produce good and feasible solutions.
Since infeasible problems are considered with an infinite objective function value, not all algorithm curves touch the upper bound 1.
This is the case for the regularized problems with $x^0 = 0$, $y^0 = e$ for all problem sizes. For the largest problems with 400 assets, even \texttt{GUROBI} with 60s limit and Kanzow--Schwartz regularization starting from $x^0 = 0$, $y^0 = 0$ were not able to reach the upper bound 1.

\begin{table}
   \caption{Results for a problem with 400 assets (pard400-e-400)}\label{tab:1}
   \scriptsize
   \tabcolsep5pt
   \centering
   \begin{tabular}{l|lll|lll|lll|lll}
      & \multicolumn{3}{c|}{VaR} & \multicolumn{3}{c|}{CVaR} & \multicolumn{3}{c|}{RVaR} & \multicolumn{3}{c}{RCVaR}\\
      $\beta$ & 0.9 & 0.95 & 0.99 & 0.9 & 0.95 & 0.99 & 0.9 & 0.95 & 0.99 & 0.9 & 0.95 & 0.99 \\
      \hline\hline
      \multicolumn{1}{l}{Alg.} & \multicolumn{12}{c}{Objective value}\\
      \hline
      1 & 62.29 & 42.79 & 58.87 & 49.58 & 55.28 & 64.67 & 64.81 & 52.08 & -- & 80.31 & 117.98 & 272.07 \\
      2 & 32.70 & 40.55 & 63.97 & 49.58 & 55.28 & 64.67 & 36.62 & 52.73 & 133.68 & 79.21 & 117.98 & 272.07 \\
      3 & 53.20 & 68.28 & 96.57 & 72.85 & 85.62 & 110.64 & 55.35 & 85.71 & 204.43 & 124.53 & 180.94 & 413.04 \\
      4 & 29.76 & 38.20 & 54.03 & 40.76 & 47.91 & 61.90 & 30.97 & 47.96 & 114.39 & 69.68 & 101.24 & 231.11 \\
      5 & 25.94 & 33.30 & 47.12 & 35.54 & 41.76 & 53.99 & 27.00 & 41.80 & 99.86 & 60.79 & -- & 201.92 \\
      6 & 29.76 & 38.20 & 54.03 & 40.76 & 47.91 & 61.90 & 30.97 & 47.95 & 114.34 & 69.68 & 101.21 & 231.02 \\
      7 & 27.30 & 35.05 & 49.58 & 37.40 & 43.94 & 56.80 & 28.41 & 44.00 & 104.92 & 63.94 & 92.86 & 201.45 \\
      8 & 29.76 & 38.20 & 54.03 & 40.76 & 47.91 & 61.90 & 30.97 & 47.95 & 114.39 & 69.68 & 101.24 & 231.11 \\
      \hline\hline
      \multicolumn{1}{l}{Alg.} & \multicolumn{12}{c}{Relative gap}\\
      \hline
      1 & 1.40 & 0.29 & 0.25 & 0.39 & 0.32 & 0.20 & 1.40 & 0.25 & -- & 0.32 & 0.27 & 0.35 \\
      2 & 0.26 & 0.22 & 0.36 & 0.39 & 0.32 & 0.20 & 0.36 & 0.26 & 0.34 & 0.30 & 0.27 & 0.35 \\
      3 & 1.05 & 1.05 & 1.05 & 1.05 & 1.05 & 1.05 & 1.05 & 1.05 & 1.05 & 1.05 & 0.95 & 1.05 \\
      4 & 0.15 & 0.15 & 0.15 & 0.15 & 0.15 & 0.15 & 0.15 & 0.15 & 0.15 & 0.15 & 0.09 & 0.15 \\
      5 & 0.00 & 0.00 & 0.00 & 0.00 & 0.00 & 0.00 & 0.00 & 0.00 & 0.00 & 0.00 & -- & 0.00 \\
      6 & 0.15 & 0.15 & 0.15 & 0.15 & 0.15 & 0.15 & 0.15 & 0.15 & 0.14 & 0.15 & 0.09 & 0.15 \\
      7 & 0.05 & 0.05 & 0.05 & 0.05 & 0.05 & 0.05 & 0.05 & 0.05 & 0.05 & 0.05 & 0.00 & 0.00 \\
      8 & 0.15 & 0.15 & 0.15 & 0.15 & 0.15 & 0.15 & 0.15 & 0.15 & 0.15 & 0.15 & 0.09 & 0.15 \\
      \hline\hline
      \multicolumn{1}{l}{Alg.} & \multicolumn{12}{c}{Computation time (s)} \\
      \hline
      1 & 60 & 61 & 62 & 60 & 61 & 69 & 67 & 73 & -- & 60 & 68 & 68 \\
      2 & 300 & 300 & 258 & 300 & 300 & 300 & 300 & 300 & 189 & 288 & 302 & 300  \\
      3 & 0.03 & 0.03 & 0.03 & 0.03 & 0.03 & 0.03 & 0.03 & 0.03 & 0.02 & 0.02 & 0.03 & 0.03 \\
      4 & 0.04 & 0.04 & 0.05 & 0.05 & 0.05 & 0.05 & 0.05 & 0.04 & 0.04 & 0.05 & 0.05 & 0.05  \\
      5 & 1.03 & 0.78 & 1.09 & 1.07 & 1.06 & 1.16 & 0.97 & 1.17 & 1.14 & 0.71 & -- & 1.04  \\
      6 & 1.53 & 1.47 & 1.52 & 1.64 & 1.51 & 1.39 & 1.50 & 1.52 & 1.37 & 1.41 & 1.25 & 1.20  \\
      7 & 0.77 & 0.78 & 0.76 & 0.76 & 0.72 & 0.81 & 0.73 & 0.83 & 0.69 & 0.77 & 0.71 & 0.69  \\
      8 & 0.87 & 0.91 & 0.80 & 0.81 & 0.92 & 0.80 & 0.81 & 0.88 & 0.83 & 0.86 & 0.83 & 0.83
   \end{tabular}
\end{table}

\begin{table}
   \caption{Results for 30 instances with 200 assets}\label{tab:3}
   \scriptsize
   \tabcolsep5pt
   \centering
   \begin{tabular}{l|lll|lll|lll|lll}
      & \multicolumn{3}{c|}{VaR} & \multicolumn{3}{c|}{CVaR} & \multicolumn{3}{c|}{RVaR} & \multicolumn{3}{c}{RCVaR} \\
      $\beta$ & 0.9 & 0.95 & 0.99 & 0.9 & 0.95 & 0.99 & 0.9 & 0.95 & 0.99 & 0.9 & 0.95 & 0.99 \\
      \hline\hline
      \multicolumn{1}{l}{Alg.} & \multicolumn{12}{c}{Average relative gap} \\
      \hline
      1 & 0.126 & 0.157 & 0.143 & 0.158 & 0.171 & 0.180 & 0.164 & 0.170 & 0.154 & 0.167 & 0.179 & 0.188  \\
      2 & 0.158 & 0.149 & 0.146 & 0.133 & 0.130 & 0.155 & 0.123 & 0.139 & 0.145 & 0.143 & 0.147 & 0.164  \\
      3 & 1.091 & 1.087 & 1.093 & 1.092 & 1.088 & 1.089 & 1.073 & 1.085 & 1.096 & 1.094 & 1.092 & 1.096  \\
      4 & 0.166 & 0.163 & 0.167 & 0.166 & 0.164 & 0.165 & 0.156 & 0.162 & 0.169 & 0.167 & 0.167 & 0.169  \\
      5 & 0.000 & 0.000 & 0.000 & 0.000 & 0.000 & 0.000 & 0.000 & 0.000 & 0.000 & 0.000 & 0.000 & 0.000  \\
      6 & 0.166 & 0.163 & 0.167 & 0.166 & 0.164 & 0.165 & 0.156 & 0.162 & 0.169 & 0.167 & 0.167 & 0.169  \\
      7 & 0.004 & 0.010 & 0.016 & 0.013 & 0.020 & 0.014 & 0.008 & 0.020 & 0.017 & 0.014 & 0.015 & 0.000  \\
      8 & 0.166 & 0.163 & 0.167 & 0.166 & 0.164 & 0.165 & 0.156 & 0.162 & 0.169 & 0.167 & 0.167 & 0.169  \\
      \hline\hline
      \multicolumn{1}{l}{Alg.} & \multicolumn{12}{c}{Average computation time (s)} \\
      \hline
      1 & 60 & 60 & 60 & 60 & 60 & 60 & 67 & 67 & 67 & 60 & 60 & 60  \\
      2 & 59 & 72 & 51 & 63 & 63 & 80 & 56 & 59 & 71 & 56 & 66 & 71  \\
      3 & 0.02 & 0.02 & 0.01 & 0.01 & 0.01 & 0.01 & 0.02 & 0.02 & 0.01 & 0.01 & 0.02 & 0.02  \\
      4 & 0.02 & 0.02 & 0.02 & 0.02 & 0.02 & 0.02 & 0.02 & 0.02 & 0.02 & 0.02 & 0.02 & 0.02  \\
      5 & 0.25 & 0.26 & 0.27 & 0.26 & 0.28 & 0.25 & 0.27 & 0.28 & 0.26 & 0.27 & 0.27 & 0.26  \\
      6 & 0.31 & 0.31 & 0.34 & 0.31 & 0.31 & 0.30 & 0.31 & 0.32 & 0.33 & 0.31 & 0.32 & 0.31  \\
      7 & 0.23 & 0.20 & 0.18 & 0.18 & 0.19 & 0.19 & 0.19 & 0.20 & 0.18 & 0.19 & 0.17 & 0.19  \\
      8 & 0.21 & 0.21 & 0.22 & 0.22 & 0.21 & 0.20 & 0.21 & 0.21 & 0.20 & 0.21 & 0.20 & 0.20  \\
      \hline\hline
      \multicolumn{1}{l}{Alg.} & \multicolumn{12}{c}{Best solution found (out of 30)} \\
      \hline
      1 & 0 & 0 & 0 & 0 & 0 & 0 & 0 & 0 & 0 & 0 & 0 & 0  \\
      2 & 0 & 0 & 0 & 0 & 1 & 0 & 1 & 1 & 0 & 0 & 0 & 0  \\
      3 & 0 & 0 & 0 & 0 & 0 & 0 & 0 & 0 & 0 & 0 & 0 & 0  \\
      4 & 0 & 0 & 0 & 0 & 0 & 0 & 0 & 0 & 0 & 0 & 0 & 0  \\
      5 & 29 & 29 & 28 & 28 & 28 & 28 & 26 & 28 & 27 & 28 & 27 & 30  \\
      6 & 0 & 0 & 0 & 0 & 0 & 0 & 1 & 0 & 0 & 0 & 0 & 0  \\
      7 & 17 & 20 & 14 & 19 & 16 & 14 & 16 & 15 & 11 & 14 & 11 & 9  \\
      8 & 0 & 0 & 0 & 0 & 0 & 0 & 0 & 0 & 0 & 0 & 0 & 0  \\
      \hline\hline
      \multicolumn{1}{l}{Alg.} & \multicolumn{12}{c}{Solution was infeasible (out of 30)}\\
      \hline
      1 & 0 & 0 & 0 & 0 & 0 & 0 & 0 & 0 & 0 & 0 & 0 & 0  \\
      2 & 0 & 0 & 0 & 0 & 0 & 0 & 0 & 0 & 0 & 0 & 0 & 0  \\
      3 & 0 & 0 & 0 & 0 & 0 & 0 & 0 & 0 & 0 & 0 & 0 & 0  \\
      4 & 0 & 0 & 0 & 0 & 0 & 0 & 0 & 0 & 0 & 0 & 0 & 0  \\
      5 & 1 & 1 & 1 & 1 & 1 & 0 & 3 & 1 & 0 & 0 & 0 & 0  \\
      6 & 0 & 0 & 0 & 0 & 0 & 0 & 0 & 0 & 0 & 0 & 0 & 0  \\
      7 & 6 & 4 & 5 & 3 & 5 & 5 & 5 & 5 & 13 & 6 & 14 & 21  \\
      8 & 0 & 0 & 0 & 0 & 0 & 0 & 0 & 0 & 0 & 0 & 0 & 0
   \end{tabular}
\end{table}

\begin{table}
   \caption{Results for 30 instances with 300 assets}\label{tab:4}
   \scriptsize
   \tabcolsep5pt
   \centering
   \begin{tabular}{l|lll|lll|lll|lll}
      & \multicolumn{3}{c|}{VaR} & \multicolumn{3}{c|}{CVaR} & \multicolumn{3}{c|}{RVaR} & \multicolumn{3}{c}{RCVaR} \\
      $\beta$ & 0.9 & 0.95 & 0.99 & 0.9 & 0.95 & 0.99 & 0.9 & 0.95 & 0.99 & 0.9 & 0.95 & 0.99  \\
      \hline\hline
      \multicolumn{1}{l}{Alg.} & \multicolumn{12}{c}{Average relative gap}\\\hline
      1 & 0.256 & 0.245 & 0.224 & 0.219 & 0.232 & 0.202 & 0.204 & 0.196 & 0.217 & 0.220 & 0.248 & 0.232  \\
      2 & 0.209 & 0.224 & 0.216 & 0.206 & 0.197 & 0.187 & 0.229 & 0.200 & 0.215 & 0.207 & 0.234 & 0.230  \\
      3 & 1.093 & 1.086 & 1.082 & 1.083 & 1.082 & 1.085 & 1.093 & 1.082 & 1.094 & 1.090 & 1.094 & 1.094  \\
      4 & 0.170 & 0.166 & 0.163 & 0.164 & 0.164 & 0.165 & 0.170 & 0.164 & 0.171 & 0.168 & 0.171 & 0.171  \\
      5 & 0.001 & 0.000 & 0.000 & 0.000 & 0.000 & 0.000 & 0.001 & 0.000 & 0.000 & 0.000 & 0.000 & 0.000  \\
      6 & 0.170 & 0.166 & 0.163 & 0.164 & 0.164 & 0.165 & 0.170 & 0.164 & 0.171 & 0.168 & 0.169 & 0.171  \\
      7 & 0.016 & 0.017 & 0.011 & 0.021 & 0.014 & 0.014 & 0.018 & 0.014 & 0.018 & 0.020 & 0.022 & 0.008  \\
      8 & 0.170 & 0.166 & 0.163 & 0.164 & 0.164 & 0.165 & 0.170 & 0.164 & 0.171 & 0.168 & 0.171 & 0.168  \\
      \hline\hline
      \multicolumn{1}{l}{Alg.} & \multicolumn{12}{c}{Average computation time (s)} \\
      \hline
      1 & 60 & 60 & 61 & 60 & 60 & 60 & 67 & 67 & 64 & 60 & 60 & 60  \\
      2 & 141 & 135 & 132 & 148 & 158 & 127 & 135 & 149 & 139 & 144 & 125 & 156  \\
      3 & 0.02 & 0.02 & 0.02 & 0.02 & 0.02 & 0.02 & 0.02 & 0.02 & 0.02 & 0.02 & 0.03 & 0.02  \\
      4 & 0.03 & 0.03 & 0.03 & 0.03 & 0.03 & 0.03 & 0.03 & 0.03 & 0.03 & 0.03 & 0.03 & 0.03  \\
      5 & 0.49 & 0.48 & 0.56 & 0.53 & 0.58 & 0.53 & 0.49 & 0.56 & 0.56 & 0.59 & 0.57 & 0.55  \\
      6 & 0.69 & 0.68 & 0.69 & 0.68 & 0.67 & 0.69 & 0.70 & 0.70 & 0.73 & 0.72 & 0.69 & 0.67  \\
      7 & 0.47 & 0.43 & 0.40 & 0.40 & 0.41 & 0.38 & 0.42 & 0.42 & 0.39 & 0.43 & 0.39 & 0.41  \\
      8 & 0.47 & 0.46 & 0.44 & 0.44 & 0.45 & 0.45 & 0.46 & 0.47 & 0.43 & 0.44 & 0.44 & 0.44  \\
      \hline\hline
      \multicolumn{1}{l}{Alg.} & \multicolumn{12}{c}{Best solution found (out of 30)} \\
      \hline
      1 & 0 & 0 & 0 & 0 & 0 & 0 & 0 & 0 & 0 & 0 & 0 & 0  \\
      2 & 0 & 0 & 0 & 0 & 0 & 0 & 0 & 0 & 0 & 0 & 0 & 0  \\
      3 & 0 & 0 & 0 & 0 & 0 & 0 & 0 & 0 & 0 & 0 & 0 & 0  \\
      4 & 0 & 0 & 0 & 0 & 0 & 0 & 0 & 0 & 0 & 0 & 0 & 0  \\
      5 & 29 & 28 & 21 & 28 & 26 & 23 & 29 & 25 & 25 & 25 & 26 & 22  \\
      6 & 0 & 0 & 0 & 1 & 0 & 0 & 0 & 0 & 0 & 0 & 0 & 0  \\
      7 & 13 & 10 & 18 & 8 & 9 & 13 & 14 & 8 & 11 & 11 & 10 & 11  \\
      8 & 0 & 0 & 0 & 0 & 0 & 0 & 0 & 0 & 0 & 0 & 0 & 0  \\
      \hline\hline
      \multicolumn{1}{l}{Alg.} & \multicolumn{12}{c}{Solution was infeasible (out of 30)} \\
      \hline
      1 & 0 & 0 & 0 & 0 & 0 & 0 & 0 & 0 & 0 & 0 & 0 & 0  \\
      2 & 0 & 0 & 0 & 0 & 0 & 0 & 0 & 0 & 0 & 0 & 0 & 0  \\
      3 & 0 & 0 & 0 & 0 & 0 & 0 & 0 & 0 & 0 & 0 & 0 & 0  \\
      4 & 0 & 0 & 0 & 0 & 0 & 0 & 0 & 0 & 0 & 0 & 0 & 0  \\
      5 & 0 & 2 & 2 & 2 & 2 & 2 & 0 & 3 & 0 & 1 & 0 & 0  \\
      6 & 0 & 0 & 0 & 0 & 0 & 0 & 0 & 0 & 0 & 0 & 1 & 0  \\
      7 & 1 & 4 & 6 & 4 & 3 & 11 & 1 & 3 & 13 & 12 & 12 & 17  \\
      8 & 0 & 0 & 0 & 0 & 0 & 0 & 0 & 0 & 0 & 0 & 0 & 1
   \end{tabular}
\end{table}

\begin{table}
   \caption{Results for 30 instances with 400 assets}\label{tab:5}
   \scriptsize
   \tabcolsep5pt
   \centering
   \begin{tabular}{l|lll|lll|lll|lll}
      & \multicolumn{3}{c|}{VaR} & \multicolumn{3}{c|}{CVaR} & \multicolumn{3}{c|}{RVaR} & \multicolumn{3}{c}{RCVaR} \\
      $\beta$ & 0.9 & 0.95 & 0.99 & 0.9 & 0.95 & 0.99 & 0.9 & 0.95 & 0.99 & 0.9 & 0.95 & 0.99  \\
      \hline\hline
      \multicolumn{1}{l}{Alg.} & \multicolumn{12}{c}{Average relative gap} \\
      \hline
      1 & 0.266 & 0.259 & 0.235 & 0.226 & 0.277 & 0.241 & 0.311 & 0.228 & 0.209 & 0.253 & 0.288 & 0.269  \\
      2 & 0.205 & 0.238 & 0.244 & 0.198 & 0.215 & 0.214 & 0.206 & 0.199 & 0.226 & 0.212 & 0.258 & 0.250  \\
      3 & 1.180 & 1.191 & 1.195 & 1.188 & 1.172 & 1.201 & 1.182 & 1.178 & 1.201 & 1.201 & 1.197 & 1.200  \\
      4 & 0.171 & 0.177 & 0.179 & 0.175 & 0.168 & 0.183 & 0.173 & 0.171 & 0.183 & 0.183 & 0.181 & 0.182  \\
      5 & 0.000 & 0.000 & 0.000 & 0.000 & 0.000 & 0.000 & 0.000 & 0.000 & 0.000 & 0.000 & 0.000 & 0.001  \\
      6 & 0.171 & 0.177 & 0.179 & 0.175 & 0.167 & 0.183 & 0.172 & 0.170 & 0.183 & 0.183 & 0.181 & 0.183  \\
      7 & 0.016 & 0.023 & 0.017 & 0.016 & 0.011 & 0.017 & 0.015 & 0.011 & 0.006 & 0.013 & 0.003 & 0.000  \\
      8 & 0.171 & 0.177 & 0.179 & 0.175 & 0.167 & 0.183 & 0.173 & 0.171 & 0.180 & 0.183 & 0.181 & 0.183  \\
      \hline\hline
      \multicolumn{1}{l}{Alg.} & \multicolumn{12}{c}{Average computation time (s)} \\
      \hline
      1 & 62 & 61 & 67 & 62 & 62 & 64 & 67 & 67 & 65 & 63 & 64 & 63  \\
      2 & 231 & 243 & 220 & 201 & 206 & 205 & 217 & 221 & 224 & 222 & 183 & 210 \\
      3 & 0.03 & 0.03 & 0.03 & 0.03 & 0.03 & 0.03 & 0.03 & 0.03 & 0.03 & 0.03 & 0.04 & 0.03  \\
      4 & 0.05 & 0.05 & 0.05 & 0.04 & 0.05 & 0.05 & 0.05 & 0.05 & 0.04 & 0.04 & 0.05 & 0.04  \\
      5 & 0.94 & 1.03 & 1.09 & 1.04 & 1.02 & 1.04 & 0.93 & 1.04 & 1.07 & 1.07 & 1.06 & 1.06  \\
      6 & 1.36 & 1.36 & 1.38 & 1.38 & 1.36 & 1.36 & 1.37 & 1.50 & 1.32 & 1.37 & 1.32 & 1.30  \\
      7 & 0.85 & 0.77 & 0.76 & 0.75 & 0.76 & 0.77 & 1.00 & 0.76 & 0.73 & 0.73 & 0.74 & 0.72  \\
      8 & 0.86 & 0.86 & 0.84 & 0.84 & 0.83 & 0.84 & 0.85 & 0.85 & 0.82 & 0.85 & 0.83 & 0.79  \\
      \hline\hline
      \multicolumn{1}{l}{Alg.} & \multicolumn{12}{c}{Best solution found (out of 30)} \\
      \hline
      1 & 0 & 0 & 0 & 0 & 0 & 0 & 0 & 1 & 0 & 0 & 0 & 0  \\
      2 & 0 & 0 & 1 & 0 & 0 & 0 & 0 & 1 & 0 & 0 & 0 & 0  \\
      3 & 0 & 0 & 0 & 0 & 0 & 0 & 0 & 0 & 0 & 0 & 0 & 0  \\
      4 & 0 & 0 & 0 & 0 & 0 & 0 & 0 & 0 & 0 & 0 & 0 & 0  \\
      5 & 24 & 25 & 19 & 22 & 13 & 20 & 24 & 14 & 16 & 19 & 16 & 21  \\
      6 & 0 & 0 & 0 & 0 & 0 & 0 & 0 & 0 & 0 & 0 & 0 & 0  \\
      7 & 14 & 8 & 12 & 9 & 16 & 12 & 12 & 15 & 16 & 14 & 16 & 11  \\
      8 & 0 & 0 & 0 & 0 & 1 & 0 & 0 & 0 & 0 & 0 & 0 & 0  \\
      \hline\hline
      \multicolumn{1}{l}{Alg.} & \multicolumn{12}{c}{Solution was infeasible (out of 30)} \\
      \hline
      1 & 0 & 0 & 0 & 0 & 0 & 0 & 0 & 0 & 20 & 0 & 0 & 1  \\
      2 & 0 & 0 & 0 & 0 & 0 & 0 & 0 & 0 & 1 & 0 & 0 & 1  \\
      3 & 0 & 0 & 0 & 0 & 0 & 0 & 0 & 0 & 0 & 0 & 0 & 0  \\
      4 & 0 & 0 & 0 & 0 & 0 & 0 & 0 & 0 & 0 & 0 & 0 & 0  \\
      5 & 1 & 1 & 1 & 1 & 5 & 1 & 2 & 5 & 0 & 0 & 2 & 0  \\
      6 & 0 & 0 & 0 & 0 & 0 & 0 & 1 & 0 & 0 & 0 & 0 & 0  \\
      7 & 2 & 1 & 12 & 2 & 5 & 12 & 4 & 5 & 12 & 11 & 13 & 19  \\
      8 & 1 & 1 & 2 & 1 & 1 & 0 & 1 & 1 & 2 & 0 & 2 & 0  \\
   \end{tabular}
\end{table}

\begin{figure}
   \centering
   \input{performance_portfolio_meanVar_kappa=10_n=200.tikz}
   \caption{Performance plot of the objective function values for $n = 200$ assets}\label{fig:1}
\end{figure}

\begin{figure}
   \centering
   \input{performance_portfolio_meanVar_kappa=10_n=300.tikz}
   \caption{Performance plot of the objective function values for $n = 300$ assets}\label{fig:2}
\end{figure}

\begin{figure}
   \centering
   \input{performance_portfolio_meanVar_kappa=10_n=400.tikz}
   \caption{Performance plot of the objective function values for $n = 400$ assets}\label{fig:3}
\end{figure}


\section{Conclusions}

We adapted the Scholtes regularization method to optimization problems with cardinality constraints and proved its convergence to S-stationary points, which is stronger than the corresponding result known for MPCCs.
Additionally, we verified that the corresponding regularized problems have better properties than the original one.
We discussed several possible risk measures for portfolio optimization under the assumption of normality and distributional ambiguity.
Finally, we compared several solution algorithms applied to these portfolio optimization problems and showed that the Scholtes regularization can keep up even with the commercial solver \texttt{GUROBI}, at least for fast, possibly local solutions.

Future research will be devoted to developing a global solution strategy based on several starting points and combinations of the proposed methods.

\subsubsection*{Acknowledgments}

Martin Branda and Michal \v{C}ervinka would like to acknowledge the support of the Czech Science Foundation (GA \v{C}R) under the grants GA13-01930S and GA15-00735S.
The work of Alexandra Schwartz and Max Bucher is supported by the 'Excellence Initiative' of the German Federal and State Governments and the Graduate School of Computational Engineering at Technische Universit{\"a}t Darmstadt.


\begin{thebibliography}{99}

\bibitem{AB} Adam, L.; Branda, M.: Nonlinear Chance Constrained Problems: Optimality Conditions, Regularization and Solvers. \textit{J. Optim. Theory Appl.} \textbf{170}(2), 419--436, 2016.

\bibitem{A} Artzner, P.; Delbaen, F.; Eber, J.-M.; Heath, D.: Coherent measures of risk. \textit{Mathematical Finance} \textbf{9}(3), 203--228, 1999.


\bibitem{BeE 13} Beck, A.; Eldar, Y.C.: Sparsity constrained optimization: optimality conditions and algorithms. \textit{SIAM J. Optim.} \textbf{23}, 1480--1509, 2013.


\bibitem{BS} Bertsimas, D.; Shioda, R.: Algorithm for cardinality-constrained quadratic optimization. \textit{Comput. Optim. Appl.} \textbf{43}, 1--22, 2009.

\bibitem{BT} Bertsimas, D.; Takeda, A.: Optimizing over coherent risk measures and non-convexities: a robust mixed integer optimization approach. \textit{Comput. Optim. Appl.} \textbf{62}, 613--639, 2015.

\bibitem{B} Bienstock, D.: Computational study of a family of mixed-integer quadratic programming problems. \textit{Math. Program.} \textbf{74}, 121--140, 1995.

\bibitem{BM} Borchers, B.; Mitchell, J.E.: An improved branch and bound algorithm for mixed integer nonlinear programs. \textit{Comput. Oper. Res.} \textbf{21}, 359--367, 1994.

\bibitem{BKS2} Burdakov, O.P.; Kanzow, Ch.; Schwartz, A.: On a reformulation of mathematical programs with cardinality constraints. In \textit{Advances in Global Optimization}, D.Y. Gao, N. Ruan, amd W.X. Xing (eds), Proceedings of the 3rd World Congress of Global Optimization (Huangshan, China, 2013), Springer, New York, 121--140, 2015.

\bibitem{BKS1} Burdakov, O.P.; Kanzow, Ch.; Schwartz, A.: Mathematical programs with cardinality constraints: Reformulation by complementarity-type conditions and a regularization method, \textit{SIAM J. Optim.} \textbf{26}(1), 397--425, 2016.

\bibitem{CKS}  \v{C}ervinka, M.; Kanzow, Ch.; Schwartz, A.: Constraint qualifications and optimality conditions for optimization problems with cardinality constraints. \textit{Math. Program., Ser. A} \textbf{160}(1), 353--377, 2016.

\bibitem{ChMBS} Chang, T.-J.; Meade, N.; Beasley, J.E.; Sharaiha, Y.M.: Heuristics for cardinality constrained portfolio optimization. \textit{Comput. Oper. Res.} \textbf{27}, 1271--1302, 2000.

\bibitem{ChHZ} Chen, L.; He, S.; Zhang, S.: Tight Bounds for some risk measures, with applications to robust portfolio selection, \textit{Oper. Res.} \textbf{59}(4), 847--865, 2011.

\bibitem{ChZ} Chopra, V.K.; Ziemba, W.T.: The effect of errors in means, variances and covariances on optimal portfolio choice. \textit{Journal of Portfolio Management} \textbf{19}(2), 6--11, 1993.

\bibitem{DY} Delage, E.; Ye, Y.: Distributionally robust optimization under moment uncertainty with application to data-driven problems. \textit{Oper. Res.} \textbf{58}(3), 595--612, 2010.

\bibitem{DFN 05} Demiguel, A.V.; Friedlander, M.P.; Nogales, F.J.; Scholtes, S.: A two-sided relaxation scheme for mathematical programs with equilibrium constraints. \textit{SIAM J. Optim.} \textbf{16}, 587--609, 2005.

\bibitem{dMGNU} DeMiguel, V.; Garlappi, L.; Nogales, J.; Uppal, R.: A generalized approach to portfolio optimization: Improving performance by constraining portfolio norms. \textit{Manag. Sci.} \textbf{55}(5), 798--812, 2009.

\bibitem{dMN} DeMiguel, V.; Nogales, F.J.: Portfolio selection with robust estimation. \textit{Oper. Res.} \textbf{57} (3), 560-577, 2009.

\bibitem{DLR 12} Di Lorenzo, D.; Liuzzi, G.; Rinaldi, F.; Schoen, F.; Sciandrone, M.: A concave optimization-based approach for sparse portfolio selection. \textit{Optim. Methods Softw.} \textbf{27}, 983--1000, 2012.


\bibitem{eGOO} El Ghaoui, L.; Oks, M; Oustry, F.: Worst-case value-at-risk and robust portfolio optimization: a conic programming approach. \textit{Oper. Res.} \textbf{51}(4), 543--556, 2003.

\bibitem{FHZ} Fabozzi, F.J.; Huang, D.; Zhou, G.: Robust portfolios: contributions from operations research and finance. \textit{Ann. Oper. Res.} \textbf{176}(1), 191--220, 2010.

\bibitem{FPW} Fastrich, B.; Paterlini, S.; Winkler, P.: Constructing optimal sparse portfolios using regularization methods. \textit{Comput. Manag. Sci.} \textbf{12}(3), 417--434, 2015.

\bibitem{FG} Frangioni, A.;  Gentile, C.: SDP diagonalizations and perspective cuts for a class of nonseparable MIQP, \textit{Oper. Res. Lett.} \textbf{35}, 181--185, 2007.

\bibitem{FGwp} Frangioni, A.;  Gentile, C.: Mean-variance problem with minimum buy-in constraints, data and documentation, online at http://www.di.unipi.it/optimize/Data/MV.html.

\bibitem{FMPSW 13} Feng, M.; Mitchell, J.E.; Pang, J.-S.; Shen, X.; W\"achter, A.: Complementarity formulation of $ \ell_0 $-norm optimization problems.
\textit{Technical Report}, Industrial Engineering and Management Sciences, Northwestern University, Evanston, IL, USA, September 2013.

\bibitem{GMSW} Gill, P.E.; W. Murray, W.; Saunders, M.A.; Wong, E.: SNOPT 7.5 User's Manual,
\textit{CCoM Technical Report} \textbf{15-3}, Center for Computational Mathematics, University of California, San Diego.

\bibitem{GMS} Gill, P.E.; Murray, W.; Saunders, M.A.: SNOPT: An SQP algorithm for large-scale constrained optimization, \textit{SIAM Review} \textbf{47}, 99-131, 2005.

\bibitem{GI} Goldfarb, D.; Iyengar, G.: Robust portfolio selection problems. \textit{Math. Oper. Res.} \textbf{28}(1), 1--38, 2003.

\bibitem{G} Gurobi Optimization, Inc.: Gurobi Optimizer Reference Manual, online at http://www.gurobi.com, 2016.

\bibitem{HKS 13} Hoheisel, T.; Kanzow, C.; Schwartz, A.: Theoretical and numerical comparison of relaxation methods for mathematical programs with complementarity constraints.
\textit{Math. Program.} \textbf{137}, 257--288, 2013.

\bibitem{KDB 09} Kadrani, A., Dussault, J.-P., Benchakroun, A.: A new regularization scheme for mathematical programs with complementarity constraints. \textit{SIAM J. Optim.} \textbf{20}, 78--103, 2009.

\bibitem{KS} Kanzow, Ch.; Schwartz, A.: A new regularization method for mathematical programs with complementarity constraints with strong convergence properties. \textit{SIAM J. Optim.} \textbf{23}, 770--798, 2013.

\bibitem{KS2} Kanzow, Ch.; Schwartz, A.: The Price of Inexactness: Convergence Properties of Relaxation Methods for Mathematical Programs with Complementarity Constraints Revisited, \textit{Math. Oper. Res.} \textbf{40}(2), 253--275, 2015.

\bibitem{KKF} Kim, J.; Kim, W.; Fabozzi, F.: Recent developments in robust portfolios with a worst-case approach. \textit{J. Optim. Theory Appl.} \textbf{161} (1), 103--121, 2014.

\bibitem{Levy} Levy, H. (Ed.): \textit{Stochastic Dominance: Investment Decision Making under Uncertainty}. Springer, New York, 2006.

\bibitem{LiF 05} Lin, G.H.; Fukushima, M.: A modified relaxation scheme for mathematical programs with complementarity constraints. \textit{Ann. Oper. Res.} 133, 63--84, 2005.

\bibitem{LPR 96} Luo, Z.-Q.; Pang, J.-S.; Ralph, D.: Mathematical Programs with Equilibrium Constraints. Cambridge University Press, Cambridge, New York, Melbourne, 1996.


\bibitem{M1} Markowitz, H.M.: Portfolio selection. \textit{The Journal of Finance} \textbf{7}, 77-91, 1952.


\bibitem{M3} Michaud, R.O.: The Markowitz optimization enigma: Is optimized optimal? \textit{Financial Analysts J.} \textbf{45}(1), 31--42, 1989.

\bibitem{MS} Murray, W.; Shek, H.: A local relaxation method for the cardinality constrained portfolio optimization problem. \textit{Comput. Optim. Appl.} \textbf{53}, 681--709, 2012.

\bibitem{OKZ 98} Outrata, J.~V.; Ko{\v{c}}vara, M.; Zowe, J.: \textit{Nonsmooth Approach to Optimization Problems with Equilibrium Constraints.} Kluwer Academic Publishers, Dordrecht, The Netherlands, 1998.

\bibitem{PP} Pa\c{c}, A.B.; P{\i}nar, M.\c{C}.: Robust portfolio choice with CVaR and VaR under distribution and mean return ambiguity. \textit{TOP} \textbf{22}, 875--891, 2014.

\bibitem{Pard} Pardalos, P.M.; Rodgers, G.P.: Computing aspects of a branch and bound algorithm for quadratic zero-one programming. \textit{Computing} \textbf{144}, 45--131, 1990.


\bibitem{PW} Pflug, G.; Wozabal, D.: Ambiguity in portfolio selection. \textit{Quantitative Finance}, 435--442, 2007.

\bibitem{P} Popescu, I.: Robust mean-covariance solutions for stochastic optimization. \textit{Oper. Res.} \textbf{55}(1), 98--112, 2005.

\bibitem{RU} Rockafellar, R.T.; Uryasev, S.: Optimization of conditional value-at-risk. \textit{Journal of Risk} \textbf{2}, 21--41, 2000.

\bibitem{RU2} Rockafellar, R.T.; Uryasev, S.: Conditional Value-at-Risk for General Loss Distributions. \textit{J. Bank. Financ.} \textbf{26}, 1443--1471, 2002.

\bibitem{RGS 10} Ruiz-Torrubiano, R.; Garc\'{i}a-Moratilla, S.; Su\'{a}rez, A.: Optimization problems with cardinality constraints. In Y.~Tenne and C.-K.~Goh (eds.):
{\em Computational Intelligence in Optimization.} Springer, 105--130, 2010.

\bibitem{Sch} Scholtes, S.: Convergence properties of a regularization scheme for mathematical programs with complementarity constraints, \textit{SIAM J. Optim.} \textbf{11}, 918--936, 2001.

\bibitem{S} Shapiro, A.: \textit{Topics in stochastic programming}. CORE lecture series. Universit\'e catholique de Louvain, Louvain-la-nueve, 2011.

\bibitem{SLK} Shaw, D.X.; Liu, S.; Kopman, L.: Lagrangian relaxation procedure for cardinality-constrained porfolio optimization. \textit{Optim. Methods Softw.} \textbf{23}, 411--420, 2008.

\bibitem{StU 10} Steffensen, S.; Ulbrich, M.: A new regularization scheme for mathematical programs with equilibrium constraints. \textit{SIAM J. Optim.} \textbf{20}, 2504--2539, 2010.

\bibitem{SZL 13} Sun, X.; Zheng, X.; Li, D.: Recent advances in mathematical programming with semi-continuous variables and cardinality constraint. \textit{J. Oper. Res. Soc. China} \textbf{1}, 55--77, 2013.

\bibitem{TK} T\"{u}t\"{u}nc\"{u}, R.H.; Koening, M.: Robust asset allocation. \textit{Ann. Oper. Res.} \textbf{132}, 157--187, 2004.

\bibitem{ZSLS 13} Zheng, X.; Sun, X.; Li, D.; Sun, J.: Successive convex approximations to cardinality-constrained convex programs: a piecewise-linear DC approach. \textit{Comput. Optim. Appl.} \textbf{59}(1), 379--397, 2014.

\bibitem{ZF} Zhu, S.; Fukushima, M.: Worst-case Conditional Value-at-Risk with application to robust portfolio management. \textit{Oper. Res.} \textbf{57}(5), 1155-1168, 2009.

\end{thebibliography}
\end{document}